\documentclass{article}

\usepackage{amsmath,amssymb,amsfonts,graphicx,color,fancyhdr}
\usepackage[latin1]{inputenc}

\newtheorem{thm}{Theorem}[section]
\newtheorem{prop}[thm]{Proposition}
\newtheorem{lem}[thm]{Lemma}
\newtheorem{df}[thm]{Definition}

\newtheorem{rem}[thm]{Remark}
\newtheorem{cor}[thm]{Corollary}

\def\be#1 {\begin{equation} \label{#1}}
\newcommand{\ee}{\end{equation}}
\def\dem {\noindent {\bf Proof : }}

\def\sqw{\hbox{\rlap{\leavevmode\raise.3ex\hbox{$\sqcap$}}$%
\sqcup$}}
\def\findem{\ifmmode\sqw\else{\ifhmode\unskip\fi\nobreak\hfil
\penalty50\hskip1em\null\nobreak\hfil\sqw
\parfillskip=0pt\finalhyphendemerits=0\endgraf}\fi}

\newcommand{\mb}{\medskip\noindent}
\newcommand{\gb}{\bigskip\noindent}
\newcommand{\R}{\mathbb R}

\newcommand{\N}{\mathbb N}
\newcommand{\Z}{\mathbb Z}

\newcommand{\s}{\mathbf S}

\newcommand{\T}{{\bf T}}

\newcommand{\OQ}{\mathbf Q}
\newcommand{\OS}{\mathbf S}
\newcommand{\OP}{\mathbf P}

\title{$L^p$ estimates for non smooth bilinear Littlewood-Paley square functions on $\R$.}
\date{\today}
\author{Fr\'ed\'eric Bernicot \\ frederic.bernicot@math.u-psud.fr}

\begin {document}

\maketitle

\begin{abstract}
In this work, some non smooth bilinear analogues of linear Littlewood-Paley square functions on the real line are studied. Mainly we prove boundedness-properties in Lebesgue spaces for them. Let us consider the function $\phi_{n}$ satisfying $\widehat{\phi_n}(\xi)={\bf 1}_{[n,n+1]}(\xi)$ and consider the bilinear operator $ S_n(f,g)(x):=\int f(x+y)g(x-y) \phi_n(y) dy$. These bilinear operators are closely related to the bilinear Hilbert transforms.
Then for exponents $p,q,r\in(2,\infty)$ satisfying $\frac{1}{p}+\frac{1}{q}+\frac{1}{r}=1$, we prove that
$$ \left\| \left( \sum_{n\in \Z} \left|S_n(f,g) \right|^2 \right)^{1/2} \right\|_{L^{r'}(\R)} \lesssim \|f\|_{L^p(\R)} \|g\|_{L^q(\R)}.$$
\end{abstract}

\tableofcontents

\hspace{0.5cm}

Let us first recall linear results about smooth and non smooth Littlewood-Paley square functions. \\
We denote by $\Psi$ a smooth function satisfying $\widehat{\Psi}(0)=0$. Then the main result of Littlewood-Paley theory claims that for all exponent $p\in(1,\infty)$ there is a constant $c=c_p$ such that
\be{homothetie} \forall f\in\s(\R), \qquad \left\| \left( \sum_{n\in\Z} \left[ \int_\R \Psi\left(\frac{y}{2^n}\right)f(x-y)\frac{dy}{2^n} \right]^2 \right)^{1/2} \right\|_{p} \leq c \|f\|_{p}. \ee
Here we have used ``dilations'' in the frequency space, we can also use translations (which corresponds to modulations in the physical space) and we have the following inequality (for $p\geq 2$)~:
\be{translation} \forall f\in\s(\R), \qquad \left\| \left( \sum_{n\in\Z} \left[ \int_\R e^{iny}\Psi(y) f(x-y)dy \right]^2 \right)^{1/2} \right\|_{p} \leq c \|f\|_{p}. \ee
Then people were interested in replacing the smooth cutoffs (in the convolution) by non smooth ones. The first result in this direction is due to L. Carleson in \cite{Carleson}, where the following result is proved~:

\begin{thm} \label{thm:C} Let $p\in[2,\infty)$ be an exponent and $\pi_{I}$ be the Fourier multiplier on $\R$ associated to the characteristic function ${\bf 1}_{I}$. So $\pi_I(f)$ is the restriction of the Fourier transform of $f$ to the interval $I$.
Then there is a constant $c=c_p$ such that
$$ \forall f\in\s(\R), \qquad \left\| \left( \sum_{n\in\Z} \left| \pi_{[n,n+1]} (f)(x) \right|^2 \right)^{1/2} \right\|_{p} \leq c \|f\|_{p}.$$
\end{thm}

\mb This result was then extended by J.L. Rubio de Francia in \cite{RF}~:

\begin{thm} \label{thm:RF} We keep the same notations. Let $p\in[2,\infty)$ be an exponent. There is a constant $c=c_p$ such that for all collection ${\mathcal I}=(I)_I$ of intervals satisfying
\be{e} \sum_{I\in{\mathcal I}} {\bf 1}_{I} \lesssim 1, \ee we have
$$ \forall f\in\s(\R), \qquad \left\| \left( \sum_{I\in{\mathcal I}} \left| \pi_{I} (f)(x) \right|^2 \right)^{1/2} \right\|_{p} \leq c \|f\|_{p}.$$
The constant $c$ only depends on the exponent $p$ and on the implicit constant in (\ref{e}).
\end{thm}

\mb In the two last theorems, the restriction $p\geq 2$ is necessary, due to a well-known counter-example, taking $f$ as $\widehat{f}:={\bf 1}_{[0,N]}$ with $N$ an integer tending to infinity.

\gb Now we are interested in obtaining bilinear version of these kind of inequalities. 

\mb We first refer the reader to the work of G. Diestel in \cite{Diestel}. In this work, the author study bilinear square functions associated to the following non smooth symbols. Let $a,b\in(0,1)$ two different scales, then for $ p,q,r\in(1,\infty)$ satisfying
$$\frac{1}{r}=\frac{1}{p}+\frac{1}{q}$$
there is a constant $c=c(p,q,a,b)$ such that for all functions $f,g\in\s(\R)$~:
$$ \left\| \left( \sum_{n\in\Z} \left[ \int_{\R^2} e^{ix(\xi_1+\xi_2)} \widehat{f}(\xi_1) \widehat{g}(\xi_2) {\bf 1}_{[a^{n},a^{n-1}]}(\xi_1) {\bf 1}_{[-b^{n},b^n]}(\xi_2) d\xi \right]^2 \right)^{1/2} \right\|_{r} \leq c \|f\|_{p} \|g\|_q.$$
In this result, the author has bilinearized the important role of the point $0$ in the frequency plane both in the two variables $\xi_1$ and $\xi_2$. These square functions are associated to multiscale paraproducts.

\gb For ten years, some far more singular bilinear operators appeared with singularities along whole a line in the frequency plane. Mainly the first studied operator is the bilinear Hilbert transform (see the work of M. Lacey and C. Thiele \cite{LT1,LT2,LT3,LT4}). Then people have weaken the tools and the proof in order to obtain boundedness for general operators owning a modulation symmetry (see the work of C. Muscalu, T. Tao and C. Thiele \cite{MTT2,mtt,MTT3}, the work of  J. Gilbert and A. Nahmod \cite{GN,GN2} ...). \\
According to these recent works, we are interested in the following bilinear operation~: for $I$ an interval, the linear Fourier multiplication operator $\pi_I$, defined by
$$\pi_I(f)(x):= \int_{\R} e^{ix\xi_1} \widehat{f}(\xi_1) {\bf 1}_{I}(\xi_1) d\xi_1$$
is then replaced by the bilinear operation (keeping the same notation)~:
$$ \pi_I(f,g)(x):= \int_{\R^2} e^{ix(\xi_1+\xi_2)} \widehat{f}(\xi_1)\widehat{g}(\xi_2) {\bf 1}_I(\xi_2-\xi_1) d\xi.$$
This bilinear multiplier can be written in the physical space as ~:
$$\pi_I(f,g)(x):= \int_{\R} f(x-t) g(x+t) \widehat{{\bf 1}_{I}}(t) dt.$$

\begin{rem} From the different cited works, we know that we can consider any non degenerate singular line. The singular variable $\xi_2-\xi_1$ can be replaced by $\xi_2-\tan(\theta)\xi_1$ with an angle $\theta\in(-\pi/2,\pi/2)\setminus\{1,-\pi/4\}$. We only deal with $\theta=\pi/4$ for convenience.
\end{rem}

\mb Such a bilinear multiplier is very closed to bilinear Hilbert transforms. The ``bilinear symbol'' ${\bf 1}_I(\xi_2-\xi_1)$ has two singular points (the two extremal points of the interval $I$). We can also decompose it with two modulations and operators similar to bilinear Hilbert transforms. \\
From the previously cited works, proving boundedness in Lebesgue spaces for bilinear Hilbert transforms, we know that for exponents $p,q,r$ satisfying
$$ 0<\frac{1}{r}=\frac{1}{p}+\frac{1}{q}<\frac{3}{2}, \qquad 1<p,q\leq \infty$$
there is a constant $C=C(p,q,r)$ such that for all interval $I\subset \R$ and for all functions $f,g\in\s(\R)$
$$\left\| \pi_I(f,g)\right\|_{r} \leq C \|f\|_{p} \|g\|_q.$$
We emphasize that the constant $C$ can be chosen uniformly with respect to the interval $I$.  \\
So it is natural to hope a positive result about boundedness for bilinear square functions defined with these ``bilinear non smooth cutoffs''.

\gb Due to a personal communication of L. Grafakos, it is quite easy to prove the following result~:
Let ${\mathcal I}=(I_n)_{n\in\Z}$ be a collection of intervals (assumption (\ref{e}) is not necessary), then for the same exponents $p,q,r$ there is a constant $c=c(p,q,r)$ such that for all functions $f_n,g_n\in\s(\R)$~:
\be{aze} \left\| \left( \sum_{n} \left|\pi_{I_n}(f_n,g_n)\right|^2 \right)^{1/2}\right\|_{r} \leq C \left\|\left( \sum_{n} |f_n|^2 \right)^{1/2} \right\|_{p} \left\|\left( \sum_{n} |g_n|^2 \right)^{1/2}\right\|_q. \ee
In fact using modulations and decomposing $\pi_{[a_n,b_n]}$ by 
$$ \pi_{[a_n,b_n]} = \pi_{[a_n,\infty)} - \pi_{[b_n,\infty)},$$
(\ref{aze}) is reduced to the same estimate replacing $\pi_{I_n}$ by $\pi_{[0,\infty)}$ which is a bounded operator $n$-independent. Then the result is a consequence of a work of L. Grafakos and J.M. Martell (Theorem 9.1 of \cite{GM2}).
Now we look for similar results when we have only one function $f$ or $g$.

\gb The first result concerning such bilinear estimates is due to M. Lacey in \cite{lacey2}, where he studied a bilinear version of (\ref{translation}). It is proved that for exponents $p,q\in[2,\infty]$ satisfying
$$ \frac{1}{2}=\frac{1}{p}+\frac{1}{q}$$
there is a constant $C=C(p,q)$ such that for all functions $f,g\in\s(\R)$ we have~:
\be{translation2} \forall f\in\s(\R), \qquad \left\| \left( \sum_{n\in\Z} \left[\int_{\R} f(x-y) g(x+y) e^{iny}\Psi(y) dy \right]^2 \right)^{1/2} \right\|_{2} \leq c \|f\|_{p}\|g\|_q. \ee
This estimate can be written in the frequency space by
\be{translation3} \forall f\in\s(\R), \qquad \left\| \left( \sum_{n\in\Z} \left[\int_{\R^2} e^{ix(\xi_1+\xi_2)} \widehat{f}(\xi_1) \widehat{g}(\xi_2) \widehat{\Psi}(\xi_2-\xi_1-n) d\xi \right]^2 \right)^{1/2} \right\|_{2} \leq c \|f\|_{p}\|g\|_q. \ee
These two last estimates exactly correspond to a bilinear version of (\ref{translation}).

\mb Now using a ``linearization argument'' and boundedness for operators related to bilinear Hilbert transforms, we have the following bilinear version of (\ref{homothetie}). For exponents $p,q,r$ satisfying
$$ 0<\frac{1}{r}=\frac{1}{p}+\frac{1}{q}<\frac{2}{3}, \qquad 1<p,q\leq \infty$$
there is a constant $C=C(p,q,r)$ such that for all interval $I\subset \R$ and for all functions $f,g\in\s(\R)$
\be{homothetie2} \left\| \left( \sum_{n\in\Z} \left[ \int_{\R} f(x-y) g(x+y) \Psi\left(\frac{y}{2^n} \right) \frac{dy}{2^n} \right]^2 \right)^{1/2} \right\|_{r} \leq c \|f\|_{p}\|g\|_{q}. \ee
We have the following frequency representation of this estimate~:
\be{homothetie3} \left\| \left( \sum_{n\in\Z} \left[ \int_{\R^2} e^{ix(\xi_1+\xi_2)} \widehat{f}(\xi_1) \widehat{g}(\xi_2) \widehat{\Psi}(2^{n}(\xi_2-\xi_1)) d\xi \right]^2 \right)^{1/2} \right\|_{r} \leq c \|f\|_{p}\|g\|_{q}. \ee
These two previous results concern ``smooth'' bilinear square functions. Our aim is to obtain ``non smooth'' results corresponding to a bilinear version of Theorem \ref{thm:C} and \ref{thm:RF}.

\gb So let us introduce some notations and then we describe our main results. \\
 Let $([a_n,b_n])_{n\in\Z}$ (with $a_n<b_n \leq a_{n+1}$) be a collection of disjoint intervals.
The disjointness property is not very important. Since the use of the ``well distributed'' notion, due to Rubio de Francia in \cite{RF}, we know that we can reduce the study of square functions associated to a collection ${\mathcal I}$ satisyfing the property of bounded covering (\ref{e}) to the study of square functions associated to a ``well distributed'' collection (which is in particular a collection of disjoint intervals).

\mb The question is the following~:

\mb {\bf Question :} For which exponents $p,q,r$, is there a constant $c=c(p,q,r)$ such that for all functions $f,g\in \s(\R)$ and all sequences of disjoints intervals $([a_n,b_n])_n$, we have~:
$$ \left\| \left( \sum_{n\in \Z} \left|\pi_{[a_n,b_n]}(f,g) \right|^2 \right)^{1/2} \right\|_{r} \leq C \|f\|_p \|g\|_q \ ?$$

\mb Now we come to our main results. We give some positive results about this question but we do not completely answer it. We will prove the following bilinear version of Theorem \ref{thm:C}

\begin{thm}  \label{thm:principal} Let $2<p,q,r'<\infty$ be exponents satisfying
 $$\frac{1}{r}=\frac{1}{p}+\frac{1}{q}.$$
We assume that the sequences $(a_n)_n$ and $(b_n)$ satisfy that for all $n\in\Z$
$$ (b_n-a_n)=(b_{n-1}-a_{n-1}) \qquad (a_{n+1}-b_n)=(a_n-{b_{n-1}}).$$
Then, there is a constant $C=C(p,q,r)$ such that for all functions $f,g\in\s(\R)$
$$ \left\| \left( \sum_{n\in \Z} \left|\pi_{[a_n,b_n]}(f,g) \right|^2 \right)^{1/2} \right\|_{r} \leq C \|f\|_p \|g\|_q.$$
\end{thm}

\mb The proof is a mixture of the well-known time-frequency analysis used for the bilinear Hilbert transform and vector-valued arguments. We have also to introduce a vector-valued version of the different tools (trees, size, energy ...).

\begin{rem} Our assumptions are satisfied for the collection of intervals $([n,n+1])_n$ taking $a_n=b_{n-1}=n$. So we obtain the bilinear version of Carleson's result (Theorem \ref{thm:C}).
\end{rem}

\begin{rem} In addition, as for the linear theorem, the restriction $2\leq p,q$ is necessary. We describe in Subsection \ref{counterexample} a counter-example (which is a bilinear version of the linear counter-example). However we do not know if the assumption $2\leq r'$ is necessary to obtain continuities of the bilinear square function. It is required for our proof but we do not have arguments proving its necessity.
\end{rem}

\begin{rem} Such a result is interesting and permit us to expect new results about bilinear operators. Today we do not know how to extend our results in the case of an arbitrary collection of disjoint strips (with not necessary the same lenghts). However such a result could permit us to obtain a sufficient condition of regularity on a symbol $m$ to obtain boundedness in Lebesgue spaces for the bilinear multiplier~: 
$$(f,g) \to T_m(f,g)(x):=\int e^{ix(\alpha+\beta)} \widehat{f}(\alpha)\widehat{g}(\beta) m(\beta-\alpha) d\alpha d\beta$$
The previously cited papers deal with symbols $m$ satisfying
$$ \left|m^{(i)}(\xi)\right| \lesssim |\xi|^{-i}$$
for all integer $i\in{0,..,N}$ with a sufficiently large integer $N$. \\
The results about our new bilinear square functions could permit us to guarantee the boundedness of $T_m$ under only the assumption that $m$ has a ``bounded variation''. We refer the reader to \cite{CRS} and \cite{Lac} for similar arguments in the linear case.
\end{rem}

\mb We will prove an other square estimate.

\begin{thm}  \label{thm:principal100} Let $1<p,q\leq\infty$ be exponents satisfying
 $$0<\frac{1}{r}=\frac{1}{p}+\frac{1}{q}<\frac{3}{2}.$$
We have no assumption on our intervals $[a_n,b_n]$ (there are only disjoints).
Then, there is a constant $C=C(p,q,r)$ such that for all functions $f_n,g\in\s(\R)$
$$ \left\| \left( \sum_{n\in \Z} \left|\int f_n(x-t) g(x+t) \widehat{{\bf 1}_{[a_n,b_n]}}(t) dt \right|^2 \right)^{1/2} \right\|_{r} \leq C  \left\|\left( \sum_n |f_n|^2 \right)^{1/2}\right\|_p\|g\|_q.$$
\end{thm}

\mb Using the ``well distributed'' notion (see \cite{RF}), we obtain the following corollary~:

\begin{cor} Let ${\mathcal I}=(I)_I$ be a collection of intervals satisfying (\ref{e}). Then for $p,q,r$ exponents of Theorem \ref{thm:principal100}, there is a constant $c=c(p,q,r)$ such that for all functions $f_n,g\in\s(\R)$
$$ \left\| \left( \sum_{n\in \Z} \left|\pi_{I}(f_n,g) \right|^2 \right)^{1/2} \right\|_{r} \leq C \left\|\left( \sum_n |f_n|^2 \right)^{1/2}\right\|_p\|g\|_q .$$
\end{cor}

\gb The plan of this paper is as follows. We dedicate the two following sections to the proof of Theorem \ref{thm:principal}. In Section \ref{section1}, we explain how we reduce the desired result to the study of ``weak type''estimates for combinatorial model sums. Then in Section \ref{section2}, we explain the mixture between the classical ``time-frequency'' analysis and $l^2$-valued arguments to conclude the proof of Theorem \ref{thm:principal}.
We conclude in Subsection \ref{counterexample} by showing with a counter-example the necessity of the restriction $2\leq p,q$ in Theorem \ref{thm:principal}.
In Section \ref{section3} we study different square functions. We begin first by explaining the main modifications to prove Theorem \ref{thm:principal100}. Then we introduce other square functions, which the study is similar.

\section{Reduction to a study of combinatorial model sums.} \label{section1}

\mb  In order to prove Theorem \ref{thm:principal}, we use the ``classical'' time-frequency analysis used for this kind of bilinear
operators. 

\mb We define the singular set $\Omega$ in the frequency plane
$$ \Omega:= \bigcup_{n\in\Z} \left\{ (\xi_1,\xi_2),\ \xi_2-\xi_1=a_n \right\}\bigcup_{n\in\Z} \left\{ (\xi_1,\xi_2),\ \xi_2-\xi_1=b_n \right\}.$$
We denote $\tilde{\Omega}$ for its embedding in $\R^3$ defined by
$$ \tilde{\Omega}:= \left\{ \xi\in\R^3,\ \xi_1+\xi_2+\xi_3=0,\ (\xi_1,\xi_2)\in\Omega \right\}.$$

\mb We recall what are {\it tiles} and {\it tri-tiles} (see for example \cite{MTT3})~:

\begin{df} A {\it tile} is a rectangle (i.e. a product of two intervals) $I\times \omega$ of area one. A {\it tri-tile} $s$ is a
rectangle $s= I_s \times \omega_s$, which contains three tiles $s_{i}=I_{s_i} \times \omega_{s_i}$ for $i=1,2,3$ such that
$$ \forall i,j\in\{1,2,3\}, \qquad I_{s_i}=I_s, \quad d(\omega_{s_1} \times \omega_{s_2} \times \omega_{s_3},\tilde{\Omega}) \simeq |I_s|^{-1}$$
and such that there is one (and only one) index $n\in\Z$ satisfying
 \be{Qn} \forall\, (\xi_1,\xi_2)\in\omega_{s_1} \times \omega_{s_2}, \qquad  a_n<\xi_2-\xi_1< b_n. \ee
\end{df}

\mb Here we recall that $\tilde{\Omega}$ corresponds to the singular set in $\R^3$. So for $s$ a tri-tile, we see that the cube $\omega_{s_1} \times \omega_{s_2} \times \omega_{s_3}$ is a Whitney cube of the open set $\R^3 \setminus \tilde{\Omega}$. 

\mb We remember the concept of {\it grid} and {\it collection} of tri-tiles~:

\begin{df}
A set $\{I\}_{I\in {\mathcal I}}$ of real intervals is called a {\it grid} if for all $k\in \Z$
\be{grid} \sum_{\genfrac{}{}{0pt}{}{I\in {\mathcal I}}{2^{k-1}\leq |I|\leq 2^{k+1}}} {\bf 1}_{I} \lesssim {\bf 1}_{\R}, \ee
where the implicit constant is independent of $k$ and of the grid. So a grid has the same structure than the dyadic grid. \\
Let $\OQ$ be a set of tri-tiles. It is called a {\it collection} if
\begin{itemize}
 \item $\left\{ I_s,\ s\in \OQ\right\} \textrm{  is a grid,}$
 \item ${\mathcal J}:=\left\{ \omega_s,\ s\in \OQ \right\} \bigcup_{i=1}^{3} \left\{ \omega_{s_i},\ s\in \OQ \right\} \textrm{  is a grid,} $
 \item $\omega_{s_i} \subsetneq \varpi \in {\mathcal J} \Longrightarrow \forall\, j\in\{1,2,3\},\ \omega_{s_j} \subset \varpi. $
\end{itemize}
We set then $\OQ_n$ the set of all the tiles $s\in \OQ$ satisfying (\ref{Qn}) for the index $n$.
\end{df}

\begin{figure}[h!]
\centering
\includegraphics{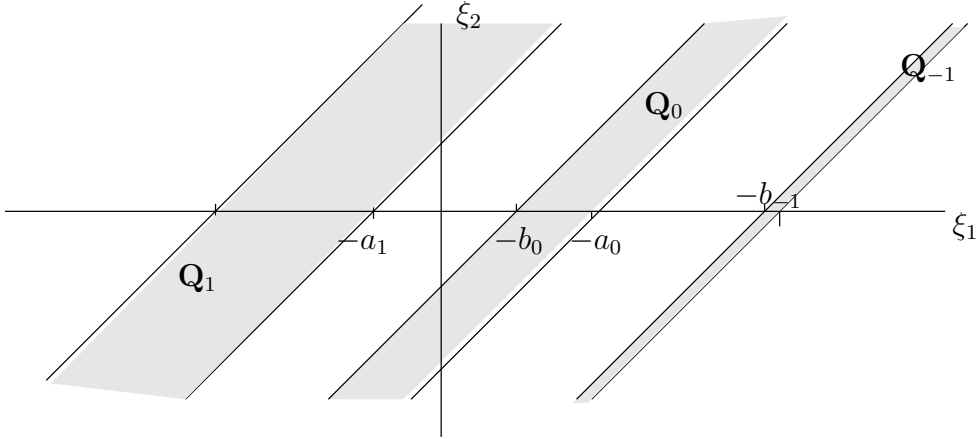}
 \caption{The strips.}
 \label{figure1}
\end{figure}

\mb We remember the notion of {\it rank one} (see for example Definition 4.9 of \cite{MTT3} for a more precise definition).

\begin{df} \label{df:rank} Let $\OQ$ be a collection of tri-tiles. It is called of {\it rank one} if 
\begin{itemize}
 \item $s\neq s'$ implies for all $j\in\{1,2,3\}$, $s_j\cap s'_j=\emptyset$.
 \item and if the collection $(I_{s'})_{s'}$ for tri-tiles $s'$ satisfying
$$ 10^7\omega_{s'_j} \supseteq \omega_{s_j} $$
is a bounded covering (with an implicit constant independent on $s$ and $j$).
\end{itemize}
\end{df}

\mb Now we can define the wave packet for a tile.

\begin{df} \label{def:wavepacket} For $P=I\times \omega$ a tile, a wave packet on $P$ is a smooth function $\Phi_P$ which has Fourier support in $\frac{9}{10} \omega$ and obeys the following estimates~: for all index $i\in\N$
$$ \left|\left[e^{ic(\omega).}\Phi_{P}\right]^{(i)}(x)\right| \lesssim |I|^{-1/2-i} \left(1+\frac{|x-c(I)|}{|I|} \right)^{-M},$$
for all exponent $M>0$ with an implicit constant depending on $M$.
For an interval $U$, we write $c(U)$ its center. So $\Phi_{P}$ is normalized in the $L^2(\R)$ space, concentrated
in space around $I$ and its spectrum is exactly contained in $\omega$.
\end{df}

\mb With classical arguments (see \cite{BG2} and \cite{BG1}), we know that we can reduce Theorem \ref{thm:principal} to the following one about model sum operators~:

\begin{thm}  \label{thm:principal2} For $2<p,q,r'<\infty$ satisfying
 $$\frac{1}{r}=\frac{1}{p}+\frac{1}{q}$$
there is a constant $C$ such that for all collection of tri-tiles $\OQ$
\begin{align*}
 \Lambda_\OQ(f,g,h)& :=\sum_{n\in\Z} \sum_{s\in \OQ_n} |I_s|^{-1/2} \left|\langle f,\Phi_{s_1}\rangle
\langle g,\Phi_{s_2}\rangle \langle h_n,\Phi_{s_3}\rangle\right| \\
 & \leq C \|f\|_p  \|g\|_q \left\| \left( \sum_{n\in Z} |h_n|^2 \right)^{1/2} \right\|_{r'}.
\end{align*}
\end{thm}

\mb As we are interested in continuities with exponents bigger than $2$, we have just to use the notion of {\it weak type}~:

\begin{df} \label{restricteddef} For $E$ a Borel set of $\R$, we write~:
$$F(E) := \left\{ f\in\s(\R),\ \forall x\in \R,\ |f(x)|\leq {\bf 1}_{E}(x) \right\}.$$ Let $p_1,p_2,p_3$ be positive exponents. We say that $\Lambda$ is {\it of weak type} $(p_1,p_2,p_3)$ if there exists
a constant $C$ such that for all measurable sets $E_1,E_2,E_3$ of finite measure with for all functions $f \in F(E_1)$, $g \in F(E_2)$ and a sequence $h:=(h_n)_n$ with $\sum_{n\in\Z} |h_n|^2 \in F(E_3)$ we have
\be{typerestreint}
\left| \Lambda(f,g,h) \right| \leq C \prod_{i=1}^{3}
|E_i|^{1/p_i}. \ee The best constant in (\ref{typerestreint}) is called
the bound of weak type and will be denoted by $C(\Lambda)$.
\end{df}

\mb By the real interpolation theory (applied to the bilinear square function) for sub-bilinear operators of weak type (see the work of L. Grafakos and N. Kalton in \cite{GK} and Exercise 1.4.17 of \cite{Gra}), Theorem \ref{thm:principal2} is reduced to the following one~:

\begin{thm} \label{thm:principal3}
Let $2< p_1,p_2,p_3<\infty$ be reals such that
$$ \frac{1}{p_1} + \frac{1}{p_2} + \frac{1}{p_3}=1$$
The trilinear form $\Lambda_\OQ$ is of weak type $(p_1,p_2,p_3)$ uniformly with respect to any finite collection $\OQ$. 
\end{thm}

\mb It is obvious that we can assume the collection $([a_n,b_n])_n$ ``well distributed'' : that is meaning
$$ \sum_{n} {\bf 1}_{\kappa[a_n,b_n]} \lesssim 1,$$ for a constant $\kappa$ as large as we want. 
In fact, due to the geometric properties of the intervals $[a_n,b_n]$, we can divide the initial collection $([a_n,b_n])_n$ with a finite number (depending on $\kappa$) of ``well distributed'' collections. So we will assume that in the following subsections.

\section{Study of these model combinatorial sums.} \label{section2}

\mb We will see that we have to use geometric properties about the strips $([a_n,b_n])_n$. Mainly we recall that we assume them to be ``well-distributed'' and we have assumed that the two length $|b_n-a_n|$ and $|a_{n+1}-b_n|$ do not depend on $n$. Let us also define a "set of reference"~:

\begin{df} We write the distance $L_1:=b_n-a_n$ and $L_2:=a_{n+1}-b_n$ (we have $L_2\gg L_1$). We will use a ``square of reference''
$$ \OS:=\left\{s\in \OQ_0,\ s_1\subset [-3L_2,3L_2]\right\}.$$
\end{df}

\mb We will use translated sets of this ``set of reference''. Let $\OQ$ be a finite collection of tri-tiles. By considering a bigger collection, we can assume the following property~: there exists an integer $N$ with
\be{propOQ} \OQ = \bigcup_{i,j=-N}^{N} \left\{ \tau_{(iL_2,jL_2,-(i+j)L_2)}(s),\  s \in \OS\cap \OQ \right\}, \ee
where for a vector $u\in\R^3$, we denote $\tau_u$ the translation operator acting on tri-tiles defined by~:
$$\tau_u(s):= I_s \times \left(u+\omega_s \right).$$

\subsection{The use of ``new trees''.}

To prove Theorem \ref{thm:principal3}, we have to organize the collection $\OQ$ with sub-collections called {\it trees}
and then study properties of orthogonality between them. The collection $\OQ$ is fixed and all the following estimates
do not depend on this collection.

\mb We recall the classical order on tiles (see Definition 4.5 in \cite{MTT3})~:

\begin{df} Let $P$ and $P'$ two tiles, we say that $P' < P$ if
$$ I_{P'} \varsubsetneq I_{P} \quad \textrm{ and } \quad 3\omega_{P} \subset 3 \omega_{P'}.$$
And we say that $P'\leq P$ if $P'<P$ or $P'=P$. We write $P' \lesssim P$ if
$$ I_{P'} \subseteq I_P \quad \textrm{ and } \quad 10^7\omega_{P} \subset 10^7 \omega_{P'}$$ and write $P' \lesssim 'P$ if $P' \lesssim P$ and $P'\nleq P$.
\end{df}

\mb We define what is a {\it tree}~:

\begin{df} Let $\T=(s_i)_i\subset \OQ$ a sub-collection of tri-tiles and $t$ an other tri-tile. $\T$ is called a $j$-tree with top $t$ (for an index $j\in\{1,2,3\}$) if there exists an index $n\in\Z$ such that $\T \subset \OQ_n$ and for all $i$~:
$$(s_i)_j \leq t_j.$$
Then we write $I_\T=I_{t}$ and $\omega_\T=\omega_t$. We say that $\T$ is a tree if it is a $j$-tree for at least one index $j\in\{1,2,3\}$.
\end{df}

\mb With this definition, it is easy to check that $\OQ$ can be divided in a collection of trees. To do this, we have just to consider the maximal tiles for the partial order.

\mb As we will see, we need to use ``vectorial trees''~:

\begin{df} Let $\T$ be a tree. For all index $k$, we will define the ``$k$-vectorized'' tree associated~:
$$ {\overrightarrow{\T}}^k := \bigcup_{s'\in\T} \left\{s\in\OQ,\ s_k=s'_k\right\}.$$
For $\T$ a tree, we define
$$ {\overrightarrow{\T}}^{1,2} := {\overrightarrow{\overrightarrow{\T}^1}}^2 = {\overrightarrow{\overrightarrow{\T}^1}}^2.$$
The second equality is due to the invariance between the different strips $D_n$~: property (\ref{propOQ}) and the fact that
$(b_n-a_n)=(a_{n+1}-b_n)=(b_{n+1}-a_{n+1})$ for all index $n\in\Z$.
\end{df}

\mb We now define several ``vectorized quantities'' (see \cite{MTT3})~:

\begin{df} Let $\OP$ be a collection of tri-tiles, $f$ a function and  $\{j,l\}=\{1,2\}$ be two indices. We define
$$ {\overrightarrow{size_j}}^l(f,\OP):= \sup_{\genfrac{}{}{0pt}{}{\T\subset \OP}{\T\ 3-\textrm{tree}}} \left( \frac{1}{|I_\T|}\sum_{s\in{\overrightarrow{\T}}^{l}} \left|\langle f,\Phi_{s_j}\rangle\right|^2 \right)^{1/2}.$$
Let $\OP$ be a collection of tri-tiles and $h:=(h_n)_{n\in\Z}$ be a sequence of  functions. We define
$${\overrightarrow{size_3}}^{1,2}(h,\OP):= \sup_{\T \subset S} \left( \frac{1}{|I_{\T}|}\sum_{n\in\Z} \sum_{s\in {\overrightarrow{\T}}^{1,2}\cap \OP_n} \left|\langle h_n,\Phi_{s_3}\rangle\right|^2 \right)^{1/2},$$
where we take the supremum over all the k-trees $\T$ with $k\neq 3$.
\end{df}

\mb We begin with important geometric remarks~:

\begin{rem} \label{rem:im2} The set $\OS$ is included in $\OQ_{0}$. We write $\tau_{(a,b,-a-b)}$ the translation in $\R^3$ of vector $(a,b,-a-b)$ in the frequency space. By the "well-distributed" property ($L_2>>L_1$), we have assumed (see (\ref{propOQ})
$$\OQ:= \bigcup_{(i,j)\in\Z} \tau_{iL_2,jL_2,-(i+j)L_2}(\OS).$$
In addition the collection $\left(\tau_{iL_2,jL_2,-(i+j)L_2}(\OS)\right)_{i,j}$ is a bounded covering of the whole collection $\OQ$ and
$$\OQ= {\overrightarrow{S}}^{1,2}.$$
For example, for $\T$ a tree we have 
$$ {\overrightarrow{\T}}^{1} = \bigcup_{j\in\Z} \tau_{0,jL_2,-jL_2}(\T).$$
\end{rem}

\begin{rem} \label{rem:imp} All the following remarks are a direct consequence of the special properties of the singular set $\Omega$~:
the strip $[a_n,b_n]$ have the same length and are separated by the same distance to the next one.
Let $\T$ be a tree then for any tree $\T'\in {\overrightarrow{\T}}^{1}$, we have ${\overrightarrow{\T}}^{1}={\overrightarrow{\T'}}^{1}$.
We have similar results for ${\overrightarrow{\T}}^{2}$ and ${\overrightarrow{\T}}^{1,2}$.
In addition there exists one and only one tree $\pi(\T) \subset {\overrightarrow{\T}}^{1,2}\cap \OS$ such that
$$ {\overrightarrow{\pi(\T)}}^{1,2}={\overrightarrow{\T}}^{1,2}.$$
We call $\pi(\T)$ the projection of $\T$ on the set $\OS$.
This projection satisfies the following property~: for all function $f_1\in\s(\R)$
$$ \overrightarrow{size_1}^2(f_1,\pi(\T)) = \overrightarrow{size_1}^2(f_1,\T)$$
and for all function $f_2\in\s(\R)~$
$$ \overrightarrow{size_2}^1(f_2,\pi(\T)) = \overrightarrow{size_2}^1(f_2,\T).$$
\end{rem}

\mb Now to cover the whole collection $\OQ$ with trees, we have to use orthogonality between them both in frequency and physical space.

\begin{df} Let $j\in\{1,2,3\}$ and $k\neq j$ be two indices. Let $D:=(\T_i)_i$ be a collection of $k$-trees with $\T_i\subset \OS$. We say that $D$ is $j$-disjoint if
\begin{enumerate}
 \item for all $i\neq i'$, for all $s\in \T_i$ and $s'\in \T_{i'}$ we have $s_j \cap s'_j= \emptyset$
 \item if for $i\neq i'$ and $s\in \T_i$, $s'\in \T_{i'}$ we have $10\omega_{s_j} \cap 10\omega_{s'_j} \neq \emptyset$ then $I_{s'} \cap I_\T = \emptyset$.
\end{enumerate}
\end{df}

\mb This definition corresponds to the ``classical one''. The restriction to trees included in $S$ may seem strange. However Remark \ref{rem:imp} explains that we can use the projection $\pi(\T)$ of a tree $\T$ on the set $\OS$ to compute a ``vectorized size''.

\begin{rem} \label{rem:disjointvect}
For $\{j,l\}:=\{1,2\}$ be indices, let $D:=(\T_i)_i$ be a $j$-disjoint collection of $k$-trees with $\T_i\subset \OS$ and $k\neq j$. Then due to the "well-distributed" property, $D$ can be decomposed in a finite number collections satisfying the following disjointness property for the vectorized trees~:
\begin{enumerate}
 \item for all $i\neq i'$, for all $s\in {\overrightarrow{\T_i}}^l$ and $s'\in {\overrightarrow{\T_{i'}}}^l$ we have $s_j \cap s'_j= \emptyset$
 \item if for $i\neq i'$ and $s\in {\overrightarrow{\T_i}}^l$, $s'\in {\overrightarrow{\T_{i'}}}^l$ we have $10\omega_{s_j} \cap 10\omega_{s'_j} \neq \emptyset$ then $I_{s'} \cap I_\T = \emptyset$.
\end{enumerate}
These properties are a direct consequence of the definition of vectorized trees using translations.
\end{rem}

\begin{rem} \label{rem:disjointvect2} We have a similar result for $j=3$. Let $D:=(\T_i)_i$ be a $3$-disjoint collection of $k$-trees with $\T_i\subset \OS$ and $k\neq 3$. Then due to the "well-distributed" property, $D$ can be decomposed in a finite number collections satisfying the following disjointness property for the vectorized trees~:
\begin{enumerate}
 \item for all $i\neq i'$ and $n\in\Z$, for all $s\in {\overrightarrow{\T_i}}^{1,2} \cap \OP_n$ and $s'\in {\overrightarrow{\T_{i'}}}^{1,2}\cap \OP_n$ we have $s_3 \cap s'_3= \emptyset$
 \item if for $i\neq i'$ and $s\in {\overrightarrow{\T_i}}^{1,2}$, $s'\in {\overrightarrow{\T_{i'}}}^{1,2}$ we have $10\omega_{s_j} \cap 10\omega_{s'_j} \neq \emptyset$ and there exists $n$ with $s_j,s_i\in\OP_n$ then $I_{s'} \cap I_\T = \emptyset$.
\end{enumerate}
\end{rem}

\mb Now we can define the other quantity ``energy'', which permits to control the $L^2$ norms of a function on a collection~:

\begin{df} Let $\OP$ be a collection, $\{j,l\}:=\{1,2\}$ be indices and $f\in\s(\R^)$ be a function. We define the ``vectorized energy''
$$ \widetilde{{\overrightarrow{energy_j}}^l}(f,\OP):=\sup_{k\in\Z} \sup_{D} 2^k \left( \sum_{\T \in D} |I_\T| \right)^{1/2},$$
where we take the supremum over all the collections $D$ of strongly $j$-disjoint trees $\T\subset \OS$ such that for all $\T\in D$
$$ \sum_{s\in {\overrightarrow{\T}}^l} \left|\langle f,\Phi_{s_j}\rangle\right|^2 \geq 2^{2k}|I_\T|$$ and for all sub-trees $\T'\subset \T$
$$ \sum_{s\in {\overrightarrow{\T'}}^l} \left|\langle f,\Phi_{s_j}\rangle\right|^2 \leq 2^{2k+2}|I_{\T'}|.$$
Similarly we define the energy for a sequence of functions $h:=(h_n)_n$ by
$$ \widetilde{{\overrightarrow{energy_3}}^{1,2}}(h,\OP):=\sup_{k\in\Z} \sup_{D} 2^k \left( \sum_{\T \in D} |I_\T| \right)^{1/2},$$
where we take the supremum over all the collections $D$ of strongly $3$-disjoint trees $\T\in \OS$ such that for all $\T\in D$
$$ \sum_{n\in\Z} \sum_{s\in {\overrightarrow{\T}}^{1,2} \cap \OQ_n} \left|\langle h_n,\Phi_{s_3}\rangle\right|^2 \geq 2^{2k}|I_\T|$$ and for all sub-trees $\T'\subset \T$
$$ \sum_{n\in \Z} \sum_{s\in {\overrightarrow{\T'}}^{1,2} \cap \OQ_n} \left|\langle h_n,\Phi_{s_3}\rangle\right|^2 \leq 2^{2k+2}|I_{\T'}|$$
\end{df}

\mb Now we want to obtain good estimates for these vectorized quantities (similarly to those obtain in the ``more classical'' case, see for example Lemmas 6.7 and 6.8 of \cite{MTT3}).

\subsection{Estimates for the vectorized ``size'' and ``energy''.}
\label{secesti}

\begin{thm} \label{thm:energy1} Let $\OP$ be a collection of tri-tiles, $f$ be a function and $\{j,l\}\in\{1,2\}$ be indices. Then we have the following estimate for the ``energy'' quantity~:
\be{energy1} \widetilde{{\overrightarrow{energy_j}}^l}(f,\OP) \lesssim \|f\|_2. \ee
\end{thm}

\dem We follow ideas of \cite{MTT3} when the authors prove their Lemma 6.7. By definition there is a collection $(\T_i)_i$ of $j$-disjoint trees of $\OS$ such that
$$ \widetilde{{\overrightarrow{energy_j}}^l}(f,\OP) \lesssim \left| \sum_{i} \sum_{s\in {\overrightarrow{\T_i}}^l} \langle f, c_s \Phi_{s_j} \rangle\right|, $$
with coefficients $c_s$ satisfying~: for all sub-trees $\widetilde{\T}\subset \T_i$
 $$ \sum_{s\in {\overrightarrow{\T'}}^l} |c_s|^2 \lesssim \frac{|I_{\T'}|}{\sum_i |I_{\T_i}|}.$$
To obtain this collection, we chose an extremizer of the definition for the vectorized energy and we chose
$$ c_s:= \frac{\langle f, \Phi_{s_j}\rangle }{2^n \left(\sum_{i} |I_{\T_i}|\right)^{1/2}}.$$
So Theorem \ref{thm:energy1} is a direct consequence of Cauchy-Schwarz inequality and the following Lemma. \findem

\begin{lem} \label{lem:ref} Let $\{j,l\}=\{1,2\}$ be indices and $D$ a collection of $j$-disjoint trees of $\OS$.
Let $c_{s_j}$ be complex numbers such that for all $\tilde{\T}$ sub-tree of $\T\in D$, we have
\be{hypp} \sum_{s\in {\overrightarrow{\tilde{\T}}}^l} |c_{s_j}|^2 \lesssim A |I_{\tilde{\T}}|. \ee
Then  we get
$$ \left\| \sum_{\T \in D} \sum_{s\in {\overrightarrow{\T}}^l} c_{s_j} \Phi_{s_j} \right\|_2^2 \lesssim A \sum_{\T\in D} |I_\T|.$$
\end{lem}

\dem We follow exactly the ideas of Lemma 6.6 of \cite{MTT3}. So we have to prove
$$ \sum_{\T,\T' \in D} \sum_{\genfrac{}{}{0pt}{}{s\in {\overrightarrow{\T}}^l}{s'\in {\overrightarrow{\tilde{T}}}^l}} \left|c_{s_j}c_{s'_j} \langle\Phi_{s_j},\Phi_{s'_j} \rangle\right|  \lesssim A \sum_{\T\in D} |I_\T|.$$
From the spectral properties of the wave packet, we sum just only the tri-tiles $s$ and $s'$ satisfying  $\omega_{s_j} \cap
\omega_{s'_j} \neq \emptyset$. By symmetry, we can assume that $|\omega_{s_j}|\leq |\omega_{s'_j}|$. From the fast spatial
decays of the wave packet, it suffices to show
\be{amontrer} \sum_{\T,\T' \in D} \sum_{\genfrac{}{}{0pt}{}{s\in
{\overrightarrow{\T}}^l}{s'\in {\overrightarrow{\T'}}^l}} \left|c_{s_j} c_{s'_j}\right| \left(\frac{|I_{s'}|}{|I_s|}
\right)^{1/2} \left( 1+ \frac{d(I_s,I_{s'})}{|I_s|} \right)^{-100}   \lesssim A \sum_{\T\in D} |I_\T|. \ee
Let us first consider the case : $|\omega_{s_j}| \simeq |\omega_{s'j}|$. We use $\left|c_{s_j}c_{s'_j}\right| \lesssim |c_{s_j}|^2 + |c_{s'_j}|^2$ and we only treat the contribution of $|c_{s_j}|^2$, the other one is similar. For a tri-tile $s$ fixed, we have seen in Remark \ref{rem:disjointvect}
that the collection $(I_{s'})_{s'}$ (for all the tri-tiles $s'$ considered in the sum) is an almost pairwise disjoint
collection and corresponds to a bounded covering. So the sum over $s'$ is bounded by a numerical constant and we also obtain the desired inequality with (\ref{hypp}). \\
Now let us consider the other case : $|\omega_{s_j}| \ll |\omega_{s'j}|$. From the assumptions, we know that for $t\in \cup_{\T\in D} \T$
$$ \sum_{s\in {\overrightarrow{t}}^l} |c_{s_j}|^2 \lesssim A|I_t|$$
and similarly for $t'$. 
As we sum under the assumption $\omega_{s_j} \subset 3 \omega_{s'_j}$, by Cauchy-Schwarz inequality, we have for $t,t'$
$$ \sum_{\genfrac{}{}{0pt}{}{s\in
{\overrightarrow{\{t\}}}^l} {\genfrac{}{}{0pt}{}{s'\in {\overrightarrow{\{t'\}}}^l}{\omega_{s_j} \subset 3 \omega_{s'_j}} }} \left|c_{s_j} c_{s'_j}\right| \lesssim \left(\sum_{s\in {\overrightarrow{t}}^l} |c_{s_j}|^2\right)^{1/2} \left(\sum_{s\in {\overrightarrow{t}}^l} |c_{s_j}|^2 \right)^{1/2} \lesssim A |I_t|^{1/2}|I_{t'}|^{1/2}.$$
So it thus suffices to show that for a fixed tree $\T\in D$
$$ \sum_{\T' \in D} \sum_{\genfrac{}{}{0pt}{}{t\in \T}{t'\in\T'}} |I_{t'}| \left( 1+ \frac{d(I_t,I_{t'})}{|I_t|} \right)^{-100}   \lesssim |I_\T|.$$
The classical observation holds : we know that there is $s\in {\overrightarrow{\{t\}}}^l$ and $s'\in {\overrightarrow{\{t'\}}}^l$ such that $\omega_{s_j} \subset 3\omega_{s'_j}$. By translation we have the same property in the set $\OS$ and so we have $\omega_{t_j} \subset 3\omega_{t'_j}$.
With classical arguments (see Lemma 6.6 of \cite{MTT3}), we can conclude that the trees $\T'$ and $\T$ are different. This is a geometric fact, we can assume ``sparseness'' (see Definition 4.4 of \cite{MTT3}) of the grid and then we have 
$$s \lesssim' t_\T \qquad s'\lesssim't_{\T'},$$
where $t_\T$ is the top of the tree $\T$. Now as $|\omega_{s_j}| \ll |\omega_{s'j}|$ and $\omega_{s_j} \subset 3\omega_{s'_j}$, we deduce that 
$\omega_{t_{\T,j}} \cap \omega_{t_{\T',j}} = \emptyset$, which permits us to prove that $\T'\neq \T$. \\
By $j$-disjointness, we obtain that the collection $(I_{t'})_{t'}$ (for all the considered tri-tiles $t'$) is a bounded covering of $I_\T^c$. We obtain also that for all $t$
\begin{align}
\sum_{\T' \in D} \sum_{t'\in\T'} |I_{t'}| \left( 1+ \frac{d(I_t,I_{t'})}{|I_t|} \right)^{-100}   & \lesssim \int_{I_\T^c} \left( 1 + \frac{d(I_t,x)}{|I_t|} \right)^{-100} dx \nonumber \\
 & \lesssim |I_t| \left( 1 + \frac{d(I_t,I_\T^c)}{|I_t|} \right)^{-10}. \label{sum}
\end{align}
By definition of a tree, a tree is always a collection of ``rank one''. Therefore for each scale $2^k \leq |I_\T|$, the collection $(I_t)_{t\in\T, |I_t|\simeq 2^k}$ is an almost disjoint collection and is a bounded covering. The desired estimate then follows by summing (\ref{sum}) for $t\in\T$. \findem

\mb Now we are interested in estimating the ``size'' quantity.

\begin{thm} \label{thm:sizevect1} Let $\OP$ be a collection of tri-tiles, $f$ be a function of $F(E_j)$ and $\{j,l\}=\{1,2\}$ be indices. Then for $p_j\geq 2$
\be{sizevect1} {\overrightarrow{size_j}}^{l}(f,\OP) \lesssim \left(\sup_{s\in \OP} \frac{1}{|I_s|} \int_{E_j}\left(1+ \frac{d(x,I_s)}{|I_s|} \right)^{-100} dx \right)^{1/p_j}. \ee
In particular we have
\be{sizevect10} \overrightarrow{size_j}^{l}(f,\OP) \lesssim \|f\|_\infty.\ee
\end{thm}

\dem We are going to apply unuseful arguments, we refer the reader to Remark \ref{rem} for a more direct proof. However the arguments, we detail, show why the restriction $p_j\geq 2$ is necessary and will be useful in Subsection \ref{subsection1}. \\
From Lemma 4.2 of \cite{MTT3b}, we know that we can compare the ``size'' quantity by introducing a vector-valued Calder\'on-Zygmund operator. We have
 \be{ineq2} \overrightarrow{size_j}^{l}(f,\OP) \simeq \sup_{\T\subset \OP}  \frac{1}{|I_\T|} \left\| \left(\sum_{s\in{\overrightarrow{\T}}^l } \left|\langle f,\Phi_{s_j}\rangle\right|^2 \frac{{\bf 1}_{I_s}}{|I_s|} \right)^{1/2} \right\|_{L^{1,\infty}(I_\T)}, \ee
where we take the supremum over all the 3-trees $\T\subset \OP$. We will follow ideas of Lemma 6.8 of \cite{MTT3}. So we fix a 3-tree $\T$ and we consider the following vector-valued operator
$$ Op_n(f):=\left(\sum_{s\in {\overrightarrow{\T}}^{l} \cap \OQ_n} \left|\langle f,\Phi_{s_j}\rangle \frac{1}{|I_s|^{1/2}}\left(1+\frac{d(x,I_s)}{|I_s|}\right)^{-200}\right|^2 \right)^{1/2}.$$
The proof of Lemma 6.8 in \cite{MTT3} explains that this vector-valued operator (for $n$ fixed) is bounded on $L^2$ and is a modulated Calder\'on-Zygmund operator as the collection ${\overrightarrow{\T}}^{l} \cap \OQ_n$ is a ``classical'' $3$-tree (we detail these claims in our Lemma \ref{lem:CZ}).
So we know that these operators $Op_n$ are $L^p$ bounded for every $p\in(1,\infty)$. In addition, obviously there is a real $\xi_0$ (which depends only on the top of the tree $\T$) such that for $s\in \T \cap \OQ_n$, the frequency interval $\omega_{s_j}$ is included into $\xi_0 + 10[a_n,b_n]:=I_n$. We write $\pi_{I_n}$ the Fourier multiplier defined by
$$ \widehat{\pi_{I_n}}(f) = {\bf 1}_{I_n} \widehat{f}$$
and so we have
$$ Op_n(f) = Op_n(\pi_{I_n}(f)).$$
From Corollary 9.4.6 of \cite{Gra}, these operators $Op_n$ (being Calder\'on-Zygmund operators) admit weighted boundedness (uniformly boundedness with respect to $n$) and so from Theorem 9.5.10 of \cite{Gra}, we have the following vector valued inequality~:
$$ \left\| \left( \sum_{n} \left|Op_n(\pi_{I_n}(f)) \right|^2 \right)^{1/2} \right\|_{p} \lesssim \left\| \left( \sum_{n} \left|\pi_{I_n}(f) \right|^2 \right)^{1/2} \right\|_{p}.$$
Then as we have assumed that the collection $([a_n,b_n])_n$ is a bounded covering, the collection $(I_n)_n$ is again a bounded covering. We can also apply the main result of R. de Francia (in \cite{RF}) about general Littlewood-Paley functionals to get for all exponent $p\geq 2$ that
$$ \left\| \left( \sum_{n} \left|\pi_{I_n}(f) \right|^2 \right)^{1/2} \right\|_{p} \lesssim \|f\|_{p}.$$
Replacing the ``smooth spatial cutoff'' by the characteristic function, we have also proved the following inequality~:
\be{ineq}  \left\| \left(\sum_{s\in{\overrightarrow{\T}}^l} \left|\langle f,\Phi_{s_j}\rangle\right|^2 \frac{{\bf 1}_{I_s}}{|I_s|} \right)^{1/2} \right\|_{p} \lesssim \|f\|_{p}. \ee
Now using (\ref{ineq2}) and H\"older inequality, we get
\begin{align*}
\frac{1}{|I_\T|} \left\| \left(\sum_{s\in{\overrightarrow{\T}}^l} \left|\langle f,\Phi_{s_j}\rangle\right|^2 \frac{{\bf 1}_{I_s}}{|I_s|} \right)^{1/2} \right\|_{L^{1,\infty}(I_\T)} & \lesssim \frac{1}{|I_\T|^{1/p}} \left\| \left(\sum_{s\in{\overrightarrow{\T}}^l} \left|\langle f,\Phi_{s_j}\rangle\right|^2 \frac{{\bf 1}_{I_s}}{|I_s|} \right)^{1/2} \right\|_{p} \\
& \lesssim \frac{1}{|I_\T|^{1/p}} \|f\|_{p}.
\end{align*}
So for $f\in F(E_j)$ supported on $I_\T$ (as it is already supported on $E_j$) we finally get
$$ \frac{1}{|I_\T|} \left\| \left(\sum_{s\in{\overrightarrow{\T}}^l} \left|\langle f,\Phi_{s_j}\rangle\right|^2 \frac{{\bf 1}_{I_s}}{|I_s|} \right)^{1/2} \right\|_{L^{1,\infty}(I_\T)} \lesssim \left( \frac{|E_j\cap I_\T|}{|I_\T|}\right)^{1/p},$$
which corresponds to what we want with $p=p_j\geq 2$.
If $f$ is supported in $2^{k+1}I_\T \setminus 2^kI_\T$ for a positive integer $k\geq 0$, then we produce the same arguments by noticing the following fact
$$\langle f,\Phi_{s_j} \rangle = 2^{-2kM} \left\langle \frac{2^{2kM}}{(1+\frac{d(.,I_\T)}{|I_\T|})^{2M}} f, (1+\frac{d(.,I_\T)}{|I_\T|})^{2M} \Phi_{s_j} \right\rangle.$$
We remark that the function $(1+\frac{d(.,I_\T)}{|I_\T|})^{2M} \Phi_{s_j}$ is always a wave packet for the tile $s_j$ (according to Definition \ref{def:wavepacket}). As $(1+\frac{d(x,I_\T)}{|I_\T|})\simeq 2^{k}$ for $x$ in the support of $f$, we produce the same arguments as before and we get the extra factor $2^{-2kM}$ with a power $M$ as large as we want. So we have proved that for every function $f$
$$\frac{1}{|I_\T|} \left\| \left(\sum_{s\in{\overrightarrow{\T}}^l} \left|\langle f,\Phi_{s_j}\rangle\right|^2 \frac{{\bf 1}_{I_s}}{|I_s|} \right)^{1/2} \right\|_{L^{1,\infty}(I_\T)} \lesssim \left(\frac{1}{|I_\T|} \int_{E_j}\left(1+ \frac{d(x,I_\T)}{|I_\T|} \right)^{-100} dx \right)^{1/p},$$
which implies the desired inequality. \findem

\mb We want to explain the study of the operator $Op_n$ when $n$ is fixed. We have used that it is a modulated vector-valued Calder\'on-Zygmund operator.

\begin{lem} \label{lem:CZ} With the notation of previous proof, let us fix an index $n$, then the operator $Op_n$ is a modulated vector-valued Calder\'on-Zygmund operator, bounded on $L^2$. All these implicit constants are uniformly bounded with respect to $n$.
\end{lem}

\dem Let us recall the $l^2$-valued operator $Op_n$~:
$$ Op_n(f):=\left(\langle f,\Phi_{s_j}\rangle \frac{1}{|I_s|^{1/2}}\left(1+\frac{d(x,I_s)}{|I_s|}\right)^{-200} \right)_{s\in{\overrightarrow{\T}}^l \cap \OQ_n}.$$
The important fact is that for $\T$ a $3$-tree, the collection ${\overrightarrow{\T}}^l \cap \OQ_n$ corresponds to a ``classical tree'' (as defined in \cite{MTT3} for example).
For the $L^2$-boundedness, we refer the reader to Exercise 10.5.7 of \cite{Gra}. We prove the $L^2$ boundedness for such an operator associated to a ``classical tree'' and so the proof holds for the collection ${\overrightarrow{\T}}^l \cap \OQ_n$.
We have also just to check the regularity about the kernel. It is given by
$$ \left|K(x,y) \right| \lesssim \left( \sum_{s\in {\overrightarrow{\T}}^l \cap \OQ_n} \left| \Phi_{s_j}(y)\right|^2 \frac{1}{|I_s|}\left(1+\frac{d(x,I_s)}{|I_s|}\right)^{-400} \right)^{1/2}.$$
We firstly use the fast spatial decay of the wave packet $\Phi_{s_j}$ and combine the different terms to get
$$ \left|K(x,y) \right| \lesssim  \left( \sum_{s\in {\overrightarrow{\T}}^l \cap \OQ_n} \frac{1}{|I_s|^2}\left(1+\frac{|x-y|}{|I_s|}\right)^{-200}\left(1+\frac{d(x,I_s)}{|I_s|}\right)^{-200} \right)^{1/2}.$$
Then we use the fact that for the scale $|I_s|\simeq 2^i$ fixed, the number of tri-tiles $s\in {\overrightarrow{\T}}^l \cap \OQ_n$ is bounded by a numerical constant. In addition the different intervals $(I_s)_s$ at a fixed length is a bounded covering of the whole space so we have~:
\begin{align*}
 \left|K(x,y) \right| & \lesssim  \left( \sum_{i\in\Z} \left(1+\frac{|x-y|}{2^{i}}\right)^{-200} \sum_{\genfrac{}{}{0pt}{}{s\in \overrightarrow{\T}^k \cap \OQ_n}{|I_s|\simeq 2^{i}}} 2^{-2i} \left(1+\frac{d(x,I_s)}{2^{i}}\right)^{-200} \right)^{1/2} \\
 & \lesssim \left( \sum_{i\in\Z} \left(1+\frac{|x-y|}{2^{i}}\right)^{-200} 2^{-2i} \right)^{1/2} \\
 & \lesssim |x-y|^{-1}.
\end{align*}
By the same estimates, we get that there exists $\xi_0$ such that for all index $a,b\in\N$
$$ \left| \partial_x^{a} \partial_y^{b} K(x,y)e^{i\xi_0y} \right| \lesssim |x-y|^{-1-a-b}, $$
which concludes the proof of the claim. The point $\xi_0$ corresponds to the top of the tree ${\overrightarrow{\T}}^l \cap \OQ_n$. \findem

\begin{rem} \label{rem} Indeed the proof of Theorem \ref{thm:sizevect1} may be simplified for the exponent $p_j=2$ due to the orthogonality properties of the space $L^2$. Then using the fact that the quantity
$$\sup_{s\in \OP} \frac{1}{|I_s|} \int_{E_j}\left(1+ \frac{d(x,I_s)}{|I_s|} \right)^{-100} dx$$ is bounded, we obtain the desired result for $p_j\geq 2$.
However this proof permits us to understand why we need the assumption $p_j\geq 2$. In the classical case (see for example \cite{MTT3}), such an estimate is proved for all exponent $p_j\geq 1$. The more complex definition of ``size'', due to the vector-valued extra term, forces us to use Littlewood-Paley inequalities for arbitrary intervals and so to assume $p_j\geq 2$. \\
This restriction is important and it explains the range for exponents in Theorem \ref{thm:principal}. This Remark will helps us in Subsection \ref{counterexample} to prove the necessity of $p_j\geq 2$ with a counter-example. \\
In addition as we will prove Theorem \ref{thm:principal} for the three exponents bigger than $2$, this precise estimate of the ``size'' quantity will be not useful. We will have just to use (\ref{sizevect10}). However to prove Theorem \ref{thm:principal100}, we will use similar arguments for exponents $p_j>1$. That is why, we prefered to explain them in details here, as we will quickly describe the arguments in Subsection \ref{subsection1} proving Theorem \ref{thm:principal100}.
\end{rem}

\mb Now we want to obtain same estimates for these quantities (``energy'' and ``size'') with the index $j=3$~:

\begin{thm} \label{thm:energy3} Let $\OP$ be a collection of tri-tiles and, $h:=(h_n)_n$ be a sequence of functions. Then we have the following estimate for the vectorized energy quantity~:
\be{energy3} \widetilde{{\overrightarrow{energy_3}}^{1,2}}(h,\OP) \lesssim \left\|\left(\sum_{n\in\Z} |h_n|^2 \right)^{1/2} \right\|_2. \ee
\end{thm}

\mb The proof is similar to the one of Theorem \ref{thm:energy1}. By the same reasoning, we have just to prove the following $l^2$-valued version of Lemma \ref{lem:ref}.

\begin{lem} \label{lem:refval} Let $D$ be a collection of $3$-disjoint trees of $\OS$.
Let $c_{s_3}$ be complex numbers such that for all $\widetilde{\T}$ sub-tree of $\T\in D$, we have
$$ \sum_{s\in {\overrightarrow{\widetilde{\T}}}^{1,2}} |c_{s_3}|^2 \lesssim A |I_{\widetilde{\T}}|.$$
Then  we get
$$ \left\| \left(\sum_{n\in\Z} \left[ \sum_{\T \in D} \sum_{s\in {\overrightarrow{\T}}^{1,2}\cap \OQ_n} c_{s_3} \Phi_{s_3} \right]^2 \right)^{1/2} \right\|_2^2 \lesssim A \sum_{\T\in D} |I_\T|.$$
\end{lem}

\dem Which is important, is that we have to control the $l^2$ norm according to the parameter $n\in\Z$ and not the $l^1$-norm.
Taking the square of the desired inequality, we have to prove
$$ \sum_{n} \sum_{\T,\T' \in D} \sum_{\genfrac{}{}{0pt}{}{s\in {\overrightarrow{\T}}^{1,2} \cap\OQ_n}{s'\in {\overrightarrow{\tilde{T}}}^{1,2} \cap \OQ_n}} \left|c_{s_3}c_{s'_3} \langle\Phi_{s_3},\Phi_{s'_3} \rangle\right|  \lesssim A \sum_{\T\in D} |I_\T|.$$
From the spectral properties of the wave packet, we sum just only the tri-tiles $s$ and $s'$ satisfying  $ \omega_{s_3} \cap
\omega_{s'_3} \neq \emptyset$. By symmetry, we may assume $|\omega_{s_3}|\leq |\omega_{s'_3}|$. From the fast spatial
decays of the wave packet, it suffices to show
\be{amontrervect} \sum_{n}\sum_{\T,\T' \in D} \sum_{\genfrac{}{}{0pt}{}{s\in
{\overrightarrow{\T}}^{1,2}\cap \OQ_n}{s'\in {\overrightarrow{\T'}}^{1,2} \cap \OQ_n}} \left|c_{s_3} c_{s'_3}\right| \left(\frac{|I_{s'}|}{|I_s|}
\right)^{1/2} \left( 1+ \frac{d(I_s,I_{s'})}{|I_s|} \right)^{-100}   \lesssim A \sum_{\T\in D} |I_\T|. \ee
Let us first consider the first case : $|\omega_{s_j}| \simeq |\omega_{s'j}|$. We use $\left|c_{s_j}c_{s'_j}\right| \lesssim |c_{s_j}|^2 + |c_{s'_j}|^2$.
We only treat the contribution of $|c_{s_j}|^2$, the other one is similar. For a tri-tile $s$ fixed, we have seen in Remark \ref{rem:disjointvect2}
that the collection $(I_{s'})_{s'}$ (for all the tri-tiles $s'$ considered in the sum) is an almost pairwise disjoint
collection and corresponds to a bounded covering, as $\OQ_n$ is a collection of rank one. So the sum over $s'$ is bounded by a numerical constant and we also obtain the desired inequality. \\
Now let us consider the other case : $|\omega_{s_j}| << |\omega_{s'j}|$. From the assumptions, we know that for $s\in
\cup \T$
$$ \sum_{t\in {\overrightarrow{\{s\}}}^{1,2}} |c_{t_j}|^2 \lesssim A|I_t|$$
and similarly for $s'$. By Cauchy-Schwarz inequality, it thus suffices to show that for a fixed tree $\T$
$$ \sum_{\T' \in D} \sum_{\genfrac{}{}{0pt}{}{s\in \T}{s'\in\T'}} |I_{s'}| \left( 1+ \frac{d(I_s,I_{s'})}{|I_s|} \right)^{-100}   \lesssim |I_\T|,$$
which was already proved in Lemma \ref{lem:ref}.
\findem

\mb We finish all these estimates by the one concerning ${\overrightarrow{size_3}}^{1,2}$.

\begin{thm} \label{thm:sizevect3}
Let $\OP$ be a collection of tri-tiles, $h:=(h_n)$ a sequence of functions satisfying  $\sum_n |h_n|^2\leq {\bf 1}_{E_3}$ for a measurable set $E_3$. Then for $p_3\geq 2$, we have
\be{sizevect3} {\overrightarrow{size_3}}^{1,2}(f,\OP) \lesssim \left(\sup_{s\in \OP} \frac{1}{|I_s|} \int_{E_3}\left(1+ \frac{d(x,I_s)}{|I_s|} \right)^{-100} dx \right)^{1/p_3}. \ee
We have above all the main estimate~:
\be{sizevect11}  {\overrightarrow{size_3}}^{1,2}(f,\OP) \lesssim \left\| \left(\sum_n |h_n|^2 \right)^{1/2} \right\|_{\infty}. \ee
\end{thm}

\dem To prove (\ref{sizevect3}), we use notations of the proof for Theorem \ref{thm:sizevect1}. With the same arguments as for Theorem \ref{thm:sizevect1}, we reduce the desired inequality to the following one~:
$$ \left\| \left( \sum_{n} \left|Op_n(h_n) \right|^2 \right)^{1/2} \right\|_{p_3} \lesssim \left\| \left( \sum_{n} \left|h_n \right|^2 \right)^{1/2} \right\|_{p_3}.$$
Now $Op_n$ is the following operator, according to a $1$ (or $2$)-tree $\T$ of $\OS$~:
$$ Op_n(h_n):=\left(\sum_{s\in {\overrightarrow{\T}}^{1,2} \cap
\OQ_n} \left|\langle h_n,\Phi_{s_3}\rangle \frac{1}{|I_s|^{1/2}}\left(1+\frac{d(x,I_s)}{|I_s|}\right)^{-200}\right|^2 \right)^{1/2}.$$
By the geometry of the strip $\OQ_n$ and of ${\overrightarrow{\T}}^{1,2}$, we can decompose $Op_n(h_n)$ with Fourier multipliers
$$Op_n(h_n):= \left(\sum_{k} \left|Op_{n,k} \pi_{k}(h_n)\right|^2 \right)^{1/2},$$
where $\pi_k$ is the Fourier multiplier by ${\bf 1}_{[kL_2,(k+1)L_2]}$ and $Op_{n,k}$ is the $l^2$-valued operator associated to the tree $\T_{n,k}$ included into 
${\overrightarrow{\T}}^{1,2} \cap\OQ_n$ and such that
$$ \forall s\in \T_{n,k}, \qquad \omega_{s_3} \subset [kL_2,(k+1)L_2].$$
The geometry of the strips implies that the tree $\T_{n,k}$ is unique and well defined by these two conditions. \\
Therefore we have now to prove
\be{ap} \left\| \left( \sum_{n,k} \left|Op_{n,k}(\pi_{k}h_n) \right|^2 \right)^{1/2} \right\|_{p_3} \lesssim \left\| \left( \sum_{n} \left|h_n \right|^2 \right)^{1/2} \right\|_{p_3}.\ee
Each operator $Op_{n,k}$ is a modulated Calder\'on-Zygmund operator (see lemma \ref{lem:CZ}), as it corresponds to a collection of tiles which is a ``classical tree''.
By using the same vector valued inequality for such operators as in Theorem \ref{thm:sizevect1}, we obtain
$$ \left\| \left( \sum_{n,k} \left|Op_{n,k}(\pi_{k}h_n) \right|^2 \right)^{1/2} \right\|_{p_3} \lesssim \left\| \left( \sum_{n,k} \left|\pi_{k} h_n \right|^2 \right)^{1/2} \right\|_{p_3}. $$
It appears the Littlewood-Paley square function, associated to the collection $([kL_2,(k+1)L_2])$, which is $L^p$-bounded only for $p\geq 2$ (see the result of R. de Francia in \cite{RF}).
So we understand as in Theorem \ref{thm:sizevect1} why we restrict the exponent $p_3$ to be bigger than $2$. Then we use Remark \ref{rem} and so we have just to prove the desired inequality for $p_3=2$. In this case (\ref{ap}) comes now easily since $([kL_2,(k+1)L_2])_k$ is a bounded covering.
Then from (\ref{ap}), the same arguments as for Theorem \ref{thm:sizevect1} hold and permit us to conclude and obtain our desired result (\ref{sizevect3}) and in particular (\ref{sizevect11}). \findem

\mb After having obtain these $l^2$-valued estimates, we apply a corresponding algorithm to conclude the proof of our main result.

\subsection{The end of the proof.} \label{subsection2}

\mb We recall the notation for $\OP$ a collection of tri-tiles~:
$$ \Lambda_\OP(f,g,h):=\sum_{n\in\Z} \sum_{s\in \OP_n} |I_s|^{-1/2} \left|\langle f,\Phi_{s_1}\rangle
\langle g,\Phi_{s_2}\rangle \langle h_n,\Phi_{s_3}\rangle\right|.$$

\gb The main result of this Section is the following one~:

\begin{thm} \label{thm:abstrait} Let $0<\theta_1,\theta_2,\theta_3$ three real numbers satisfying
$$ \theta_1+\theta_2+\theta_3=1.$$
Then there exists a constant $C=C(\theta)$ such that for all finite collections $\OQ$ of tri-tiles, all functions $f_1,f_2$ and all sequences $f_3:=(f_{3,n})_n$
\begin{align*}
\Lambda_{\OQ}(f_1,f_2,f_3) & \leq C \prod_{\{j,l\}=\{1,2\}} \left[\widetilde{{\overrightarrow{energy_j}}^l}(f_j,\OQ)\right]^{1-\theta_j} \left[{\overrightarrow{size_j}}^{l}(f_j,\OQ)\right]^{\theta_j} \\
 & \hspace{2cm} \left[\widetilde{{\overrightarrow{energy_3}}^{1,2}}(f_3,\OQ)\right]^{1-\theta_3} \left[{\overrightarrow{size_3}}^{1,2}(f_3,\OQ)\right]^{\theta_3}.
\end{align*}
\end{thm}

\mb This subsection is dedicated to prove this result. Let us first recall how we deduce Theorem \ref{thm:principal3}.

\gb {\bf Proof of Theorem \ref{thm:principal3}} This theorem is a direct consequence of Theorem \ref{thm:abstrait} (applied with $\theta_j=1-\frac{2}{p_j}\in(0,1)$). With the help of the different estimates obtained in Subsection \ref{secesti} for the vectorized ``size'' and ``energy'' quantities : (\ref{energy1}), (\ref{sizevect10}), (\ref{energy3}) and (\ref{sizevect11}), we exactly obtain what we want. \findem

\mb We have also to prove Theorem \ref{thm:abstrait}. Let us begin to estimate the trilinear form $\Lambda$ on a vectorized tree~:

\begin{prop}[Tree estimate] \label{treeesti1} Let $\T$ be a tree of a collection $\OP$ then
$$ \Lambda_{{\overrightarrow{\T}}^{1,2}}(f_1,f_2,f_3) \lesssim |I_\T| {\overrightarrow{size_1}}^2(f_1,\OP) {\overrightarrow{size_2}}^1(f_2,\OP) {\overrightarrow{size_3}}^{1,2}(f_3,\OP)  .$$
\end{prop}

\dem For example, let us assume that $\T$ is a $3$-tree.
By definition
\begin{align*}
\Lambda_{{\overrightarrow{\T}}^{1,2}}(f_1,f_2,f_3) & := \sum_{n\in\Z} \sum_{s\in {\overrightarrow{\T}}^{1,2} \cap \OQ_n} |I_s|^{-1/2} \left|\langle f_1,\Phi_{s_1}\rangle
\langle f_2,\Phi_{s_2}\rangle \langle f_{3,n},\Phi_{s_3}\rangle\right|.
\end{align*}
We write $\tau_{i,j}$ for the translation in frequency space of the vector $(i,j,-i-j)L_2 \in\R^3$. 
We have
$${\overrightarrow{\T}}^{1,2} \cap \OQ_n:= \bigcup_{j-i=c(n)} \tau_{i,j}(\T),$$
where $c(n)$ is an integer depending on $n$ and $c$ is a one-to-one map from $\Z$ to $\Z$.
Due to the geometric assumptions about the strips $\OQ_n$, we have
$\tau_{i,j}(s)=\tau_{0,j}\tau_{i,0}(s)$ and the tile $(\tau_{i,j}(s))_{1}$ is only dependent on $i$.
Also using Cauchy-Schwarz inequality, we get~:
\begin{align*}
\Lambda_{{\overrightarrow{\T}}^{1,2}}(f_1,f_2,f_3) & = \sum_{s\in \T} \sum_{i,j}  |I_s|^{-1/2} \left|\langle f_1,\Phi_{\tau_{i,0}(s)_1}\rangle
\langle f_2,\Phi_{\tau_{0,j}(s)_2}\rangle \langle f_{3,c^{-1}(i-j)},\Phi_{\tau_{i,j}(s)_3}\rangle\right| \\
 & \leq \sum_{s\in \T} \left(\sum_{i} \left|\langle f_1,\Phi_{\tau_{i,0}(s)_1}\rangle\right|^2
\right)^{1/2} \left( \sum_{j} \left|\langle f_2,\Phi_{\tau_{0,j}(s)_2}\rangle\right|^2\right)^{1/2} \\
 & \hspace{2cm} \left(\sum_{i,j} |I_s|^{-1} |\langle f_{3,c^{-1}(i-j)},\Phi_{\tau_{i,j}(s)_3}\rangle|^2 \right)^{1/2} \\
 & \leq \left(\sum_{s\in \T} \sum_{i} \left|\langle f_1,\Phi_{\tau_{i,0}(s)_1}\rangle\right|^2
\right)^{1/2} \left( \sum_{s\in \T} \sum_{j} \left|\langle f_2,\Phi_{\tau_{0,j}(s)_2}\rangle\right|^2\right)^{1/2} \\
& \hspace{2cm} \sup_{s\in \T} \left( \sum_{i,j} |I_s|^{-1} |\langle f_{3,c^{-1}(i-j)},\Phi_{\tau_{i,j}(s)_3}\rangle|^2\right)^{1/2}.
\end{align*}
Now using Remark \ref{rem:im2}, we have~:
\begin{align*}
\lefteqn{\Lambda_{{\overrightarrow{\T}}^{1,2}}(f_1,f_2,f_3) \leq \left(\sum_{s\in {\overrightarrow{\T}}^2} \left|\langle f_1,\Phi_{s_1}\rangle\right|^2 \right)^{1/2} }  & & \\
 & & \left( \sum_{s\in {\overrightarrow{\T}}^1} \left|\langle f_2,\Phi_{s_2}\rangle\right|^2\right)^{1/2}
 \sup_{s\in \T} \left( \sum_{\genfrac{}{}{0pt}{}{n\in\Z}{s'\in {\overrightarrow{ \{s\}}}^{1,2}\cap \OQ_n }} |I_s|^{-1} |\langle f_{3,n},\Phi_{s'_3}\rangle|^2\right)^{1/2}.
\end{align*}
As every singleton $\{s\}$ can be considered as a $1$ (or $2$)-tree, we finally obtain
\begin{align*}
\Lambda_{{\overrightarrow{\T}}^{1,2}}(f_1,f_2,f_3) \leq
|I_\T| {\overrightarrow{size_1}}^2(f_1,\OP) {\overrightarrow{size_2}}^1(f_2,\OP) {\overrightarrow{size_3}}^{1,2}(f_3,\OP).
\end{align*}
We have also proved the desired inequality.
The proof is exactly the same for a $1$-tree or a $2$-tree.
\findem

\mb Then we have to adapt the main combinatorial algorithm (Proposition 12.2 of \cite{MTT3}) which permit us to understand the link between the two quantities $size_j$ and $\widetilde{energy_j}$.

\mb We have the following theorem~:

\begin{thm} \label{thm:algo1} Let $\{j,l\}=\{1,2\}$ be fixed indices and $\OP$ be a collection of tiles satisfying (\ref{propOQ}) and such that
 $$ {\overrightarrow{size_j}}^{l}(f_j,\OP) \leq 2^{-n}  \widetilde{{\overrightarrow{energy_j}}^l}(f_j,\OP).$$
Then we may decompose $\OP=\OP^1 \cup \OP^2$ with a collection $\OP^1$ satisfying (\ref{propOQ}) and
 \be{prop1} {\overrightarrow{size_j}}^{l}(f_j,\OP^1) \leq 2^{-n-1}  \widetilde{{\overrightarrow{energy_j}}^l}(f_j,\OP) \ee
and $\OP^2$ is a collection of vectorial trees $\OP^2=({\overrightarrow{\T_i}}^{1,2})_i$ such that
\be{prop2} \sum_{i} |I_{\T_i}| \lesssim 2^{2n}. \ee
\end{thm}

\dem We follow ideas of Proposition 12.2 in \cite{MTT3}. Let us denote the energy $E:=\widetilde{{\overrightarrow{energy_j}}^l}(f_j,\OP)$.
We initialize $D$ a collection of trees to the empty collection.
We consider the set of all $3$ trees $\T$ of $\OP$ such that
\be{hyp1} \sum_{s\in {\overrightarrow{\T}}^l} \left|\langle f_j,\Phi_{s_j}\rangle \right|^2 \geq 4^{-1} \left(2^{-n} E \right)^{2} |I_\T| . \ee
According to Remark \ref{rem:imp}, there is a only one tree $\pi(\T) \subset \OS\cap {\overrightarrow{\T}}^{1,2}$ such that
\be{hyp2} \sum_{s\in {\overrightarrow{\pi(\T)}}^l} \left|\langle f_j,\Phi_{s_j}\rangle \right|^2 \geq 4^{-1} \left(2^{-n} E \right)^{2} |I_\T| \ee
as the two sums are equal.
Now we apply the ``classical'' algorithm to theses trees $\pi(\T)\in \OS$. We refer the reader to Proposition 12.2 in \cite{MTT3} for a precise description of the algorithm. \\
As $\OS$ is a collection of rank one (see Definition \ref{df:rank}), it works.
We choose $\pi(\T)$ among all such trees so that the center of the top $\xi_j=c(\omega_{\pi(\T),j})$ is maximal for the classical order on $\R$ and such that $\pi(\T)$ is maximal for the order on tri-tiles. We denote also $\T_1$ this tree and $\T_1'$ the following $j$-tree
$$ \T_1':=\left\{ s\in S \setminus \T_1,\ s_j\leq s_{\T_1,j} \right\}.$$
Then we add $\T_1$ and $\T_1'$ into the collection $D$ and we remove ${\overrightarrow{\T_1}}^{1,2}$ and ${\overrightarrow{\T'_1}}^{1,2}$ of the collection $\OP$. The important fact is that the new collection $\OP\setminus \left\{{\overrightarrow{\T_1}}^{1,2},{\overrightarrow{\T'_1}}^{1,2}\right\} $ satisfies (\ref{propOQ}) due to the symmetry by translation of the frequency singular space. \\
Then we repeat the algorithm and we construct two sequences of trees $(\T_i)_i$ and $(\T_i')$. When this algorithm is finished, we have construct a collection $(\T_i)_{i}$ and $(\T_i')_{i}$ of trees of $\OS$. We define $\OP_1$ and $\OP_2$ as
$$\OP_2:=\bigcup_{i} {\overrightarrow{\T_i}}^{1,2} \bigcup {\overrightarrow{\T'_i}}^{1,2}$$ and
$$ \OP_1:= \OP \setminus \OP_2.$$
By definition of the algorithm, the property (\ref{prop1}) and (\ref{propOQ}) is obvious and we have just to check (\ref{prop2}).
By classical arguments, we see that the trees $(\T_i)_{i}$ are strongly $j$-disjoints and are included in $\OS$.
In fact as the strips $[a_n,b_n]$ are supposed "well-distributed", it is quite easy to check that the vectorial tree $({\overrightarrow{\T_i}}^l)_i$ are $j$-disjoint or could be divided in a bounded number of such collection (because we have just translated the disjointness from $\OS$ to the translated sets of $\OS$, which are almost disjoint due to the "well-distributed" property of the strips, see Remark \ref{rem:disjointvect}).
So by definition of ${\overrightarrow{energy_j}}^l$ with (\ref{hyp2}), we obtain that
$$ \sum_{i} |I_{\T_i}| \lesssim 2^{2n}.$$
As the trees $\T_i$ and $\T_i'$ have the same top, we have the same property for the collection $(\T'_i)_i$, which conclude the proof.
\findem

\mb We now describe a similar result for the index $j=3$~:

\begin{thm} \label{thm:algo2} Let $\OP$ be a collection of tiles such that
 $$ {\overrightarrow{size_3}}^{1,2}(f_3,\OP) \leq 2^{-n}  \widetilde{{\overrightarrow{energy_3}}^{1,2}}(f_3,\OP),$$
where $f_3:=(f_{3,n})_n$ is a sequence of functions. \\
Then we may decompose $\OP=\OP^1 \cup \OP^2$ with
 \be{prop10} {\overrightarrow{size_3}}^{1,2}(h,\OP^1) \leq 2^{-n-1}  \widetilde{{\overrightarrow{energy_3}}^{1,2}}(h,\OP) \ee
and $\OP^2$ is a collection of vectorial trees $\OP^2=({\overrightarrow{\T_i}}^{1,2})_i$ such that
\be{prop20} \sum_{i} |I_{\T_i}| \lesssim 2^{2n}. \ee
\end{thm}

\dem The proof is similar to the previous one. With Remark \ref{rem:imp} we can deal with trees $\T$ in the set $\OS$. Then we apply the ``classical'' algorithm which works as $\OS$ is of rank one. Then (\ref{prop20}) follows from the definition of ${\overrightarrow{energy_3}}^{1,2}$ as previously. \findem

\mb We are now able to prove our main result~:

\mb {\bf Proof of Theorem \ref{thm:abstrait}} As we have seen, we can consider that $\OQ$ satisfies (\ref{propOQ}) (else we complete $\OQ$ obtaining a bigger collection satisfying this property). Then with Proposition \ref{treeesti1} and Theorems \ref{thm:algo1} and \ref{thm:algo2}, we can apply the ``classical'' arguments to obtain Theorem \ref{thm:abstrait} (see Corollary 12.3 of \cite{MTT3} for example). Fro an easy reference, we describe it. By iterating Theorems \ref{thm:algo1} and \ref{thm:algo2}, there exists a partition 
$$ \OQ = \bigcup_{n \in \Z} \OQ^n$$
where for each $n \in \Z$ and $\{j,l\}=\{1,2\}$ we have
$$ {\overrightarrow{size_j}}^{l}(f_j,\OQ^n) \leq \min\left\{2^{-n} \widetilde{{\overrightarrow{energy_j}}^l}(f_j,\OQ), {\overrightarrow{size_j}}^{l}(f_j,\OQ) \right\},$$
and for $j=3$ we have~:
$$ {\overrightarrow{size_3}}^{1,2}(f_3,\OQ^n) \leq \min\left\{2^{-n} \widetilde{{\overrightarrow{energy_3}}^{1,2}}(f_3,\OQ), {\overrightarrow{size_3}}^{1,2}(f_3,\OQ) \right\}.$$
In addition the grid $\OQ^n$ can be covered by a collection $D_n$ of trees such that
\be{arbre} \sum_{\T\in D_n} |I_\T| \lesssim 2^{2n}. \ee
With Proposition \ref{treeesti1}, we get~:
\begin{align*}
\Lambda_{\OQ^n}(f_1,f_2,f_3) & \lesssim \sum_{\T\in D_n} |I_\T| \prod_{\{j,l\}=\{1,2\}} \min\left\{2^{-n} \widetilde{{\overrightarrow{energy_j}}^l}(f_j,\OQ), {\overrightarrow{size_j}}^{l}(f_j,\OQ) \right\} \\ 
 & \hspace{2cm} \min\left\{2^{-n} \widetilde{{\overrightarrow{energy_3}}^{1,2}}(f_3,\OQ), {\overrightarrow{size_3}}^{1,2}(f_3,\OQ) \right\} \\
 & 2^{2n} \prod_{\{j,l\}=\{1,2\}} \min\left\{2^{-n} \widetilde{{\overrightarrow{energy_j}}^l}(f_j,\OQ), {\overrightarrow{size_j}}^{l}(f_j,\OQ) \right\} \\
 & \hspace{2cm} \min\left\{2^{-n} \widetilde{{\overrightarrow{energy_3}}^{1,2}}(f_3,\OQ), {\overrightarrow{size_3}}^{1,2}(f_3,\OQ) \right\}.
\end{align*}
Then we can compute the sum over $n$ and we get the desired inequality~:
\begin{align*}
\Lambda_{\OQ}(f_1,f_2,f_3) & \leq \sum_{n\in\Z} \Lambda_{\OQ^n}(f_1,f_2,f_3) \\
 & \lesssim \prod_{\{j,l\}=\{1,2\}} \widetilde{{\overrightarrow{energy_j}}^l}(f_j,\OQ)^{1-\theta_j} {\overrightarrow{size_j}}^{l}(f_j,\OQ)^{\theta_j} \\
 & \hspace{2cm} \widetilde{{\overrightarrow{energy_3}}^{1,2}}(f_3,\OQ)^{1-\theta_3} {\overrightarrow{size_3}}^{1,2}(f_3,\OQ)^{\theta_3}.
\end{align*}
 \findem

\subsection{Necessity of $p,q\geq 2$ in Theorem \ref{thm:principal}.}
\label{counterexample}

In this subsection, we want to prove the necessity of the assumption $p,q\geq 2$ in Theorem \ref{thm:principal}. About the assumption $r'\geq 2$, we do not know if it is necessary or not.

\begin{prop} In Theorem \ref{thm:principal}, the assumption $p,q\geq 2$ is necessary to have the boundedness property of the square functions.
\end{prop}

\dem As we have seen in proving Theorems \ref{thm:sizevect1} and \ref{thm:sizevect3} (see Remark \ref{rem}), the restriction for the exponents is related to the use of continuities for square Littlewood-Paley functions associated to arbitrary intervals. From the work of Rubio de Francia (in \cite{RF}), we know that this one is continuous on $L^p$ only for $p\geq 2$.
And this condition is necessary and there was a well-known example associated to the interval $([n,n+1])_{n}$. \\
As this particular collection of intervals could appear in our problem, the same counter example will work. We detail how we adapt it for our bilinear problem. \\
So assume that for $p\leq 2$ we have such an equality~:
\be{counter} \left\| \left( \sum_{n\in \Z} \left|\int f(x-t) g(x+t) \widehat{{\bf 1}_{[n,n+1]}}(t) dt \right|^2 \right)^{1/2} \right\|_{r} \leq C \|f\|_p \|g\|_q. \ee

We chose the function $g$ as $\widehat{g}={\bf 1}_{[0,1/2]}$ and $f$ such that
$ \widehat{f}={\bf 1}_{[0,2P]}$ for a large integer $P>>1$.
We denote~:
$$ T_{n}(f,g):=\int f(x-t) g(x+t) \widehat{{\bf 1}_{[n,n+1]}}(t) dt.$$
Then an easy computation gives us that for an integer $n\in[0,P]$, we have
$$ \widehat{T_n(f,g)}(\xi)=\frac{1-|\xi-n|}{2}{\bf 1}_{|\xi-n|\leq 1}.$$
Also we deduce that $T_n(f,g)(x)=e^{-inx}T_0(f,g)(x)$, hence
$$ \left\| \left( \sum_{n=0}^P \left|\int f(x-t) g(x+t) \widehat{{\bf 1}_{[n,n+1]}}(t) dt \right|^2 \right)^{1/2} \right\|_{r} =(P+1)^{1/2}\|T_0(f,g)\|_{r} \simeq P^{1/2}.$$
In the last inequality, we have used $T_0(f,g)=T_0(\widehat{{\bf1}_{[0,2]}},g)$ which is independent on $P$. \\
However by homogeneity, it comes
$$\|f\|_{p} \simeq P^{1/p'}.$$
So such an inequality as $(\ref{counter})$ would imply~:
$$ P^{1/2} \lesssim P^{1/p'}$$
for all integer $P$ as large as we want. That is why $p$ should be bigger than $2$. The same reasonning holds for $q$.
\findem

\section{Other square  functions.} \label{section3}

\subsection{Proof of Theorem \ref{thm:principal100}.} \label{subsection1}

 We will apply similar arguments. The reduction to model sums is the same and so we have to prove the ``restricted weak type estimates'' for the following trilinear forms~:
$$ \Lambda_\OQ(f,g,h):=\sum_{n\in\Z} \sum_{s\in \OQ_n} |I_s|^{-1/2} \left|\langle f_n,\Phi_{s_1}\rangle
\langle g_n,\Phi_{s_2}\rangle \langle h_n,\Phi_{s_3}\rangle\right|,$$
where $f:=(f_n)_n$ and $h:=(h_n)_n$ are twho sequences of functions.
In Theorem \ref{thm:principal100}, we have not the previous restriction for exponents. So as described in \cite{mtt,MTT3}, we need to use ``restricted weak type estimates''. 

\begin{df} \label{wrestricteddef} For $E$ a Borel set of $\R$, we write~:
$$F(E) := \left\{ f\in\s(\R),\ \forall x\in \R,\ |f(x)|\leq {\bf 1}_{E}(x) \right\}.$$ Let $p_1,p_2,p_3$ be non vanishing exponents with exactly one (and only one) negative exponent $p_\alpha$. We say that $\Lambda$ is {\it of restricted weak type} $(p_1,p_2,p_3)$ if there exists
a constant $C$ such that for all measurable sets $E_1,E_2,E_3$ of finite measure, there exists a substantial set $E'_\alpha$ (e.d.  $E'_\alpha\subset E_\alpha$ and $|E'_\alpha|\geq |E_\alpha|/2$) with for all functions $(f_n)_n,g,(h_n)_n\in\s(\R)$ with $\sum_{n} |f_n|^2 \in F(E'_1)$, $g \in F(E'_2)$ and $\sum_{n\in\Z} |h_n|^2 \in F(E'_3)$ we have
\be{wtyperestreint}
\left| \Lambda(f,g,h) \right| \leq C \prod_{i=1}^{3}
|E_i|^{1/p_i}. \ee 
We denote for $\beta\neq \alpha$, $E'_\beta=E_\beta$ for convenience. The best constant in (\ref{typerestreint}) is called
the bound of weak type and will be denoted by $C(\Lambda)$.
\end{df}

\mb According to interpolation results on ``restricted weak type estimates'' (see \cite{MTT2}), we reduce \ref{thm:principal100} to the following one~:

\begin{thm} \label{thm:principal103}
Let $\alpha\in\{1,2,3\}$ and $ p_1,p_2,p_3$ be exponents such that
$$ \frac{1}{p_1} + \frac{1}{p_2} + \frac{1}{p_3}=1, \quad  -\frac{1}{2}<\frac{1}{p_\alpha} < 0, \quad \forall \beta\neq \alpha,\ \frac{1}{2} < \frac{1}{p_\beta} < 1.$$
The trilinear form $\Lambda_\OQ$ is of weak type $(p_1,p_2,p_3)$ uniformly with respect to any finite collection $\OQ$. 
\end{thm}

\mb In this case, we do not use the ``set of reference'' $\OS$, but just use a ``strip of reference'' : $\OQ_0$.
We assume that $\OQ_0$ is the bigger strip : $L_2:=|a_0-b_0|=\max_n |a_n-b_n|$.
We remember that for this theorem, we have no geometric assumptions for the strips, they are just disjoint. \\
So we chose a collection $\OQ$ satisfying the following property
\be{propOQ100} \OQ = \bigcup_{i}^{N} \left\{ \tau_{(iL_2,0,-iL_2)}(s),\  s \in \OQ_0 \right\}, \ee
which is possible (in taking a bigger collection).

\gb We have to define (a little) different ``size'' quantities~:

\begin{df} Let $\OP$ be a collection of tri-tiles, $f:=(f_n)_n$ a sequence of functions and $j\in\{1,3\}$ and $l\in \{1,2,3\}$ be two indices (maybe the same). We define
$$ {\overrightarrow{size_j}}^l(f,\OP):= \sup_{\genfrac{}{}{0pt}{}{\T \subset \OP_0}{\T k-\textrm{tree}}} \left( \frac{1}{|I_\T|}\sum_{n} \sum_{s\in{\overrightarrow{\T}}^{l}\cap \OP_n} \left|\langle f_n,\Phi_{s_j}\rangle\right|^2 \right)^{1/2},$$
where we take the supremum over all the $k$-tree with $k\neq j$.
For $g$ a function, we similarly define
$$ {size_2}(g,\OP):= \sup_{\genfrac{}{}{0pt}{}{\T \subset \OP_0}{\T k-\textrm{tree}}} \left( \frac{1}{|I_\T|} \sum_{s\in \T} \left|\langle g,\Phi_{s_2}\rangle\right|^2 \right)^{1/2},$$
where we take the supremum over all the $k$-tree with $k\neq 2$.
\end{df}

\mb We define ``energy'' quantities~:

\begin{df} Let $\OP$ be a collection, $f:=(f_n)_n$ be a sequence of functions and $j\in\{1,3\}$ and $l\in \{1,2,3\}$ be two indices. We define the ``energy''
$$ \widetilde{{\overrightarrow{energy_j}}^l}(f,\OP):=\sup_{k\in\Z} \sup_{D} 2^k \left( \sum_{\T \in D} |I_\T| \right)^{1/2},$$
where we take the supremum over all the collections $D$ of strongly $j$ disjoint trees $\T\subset \OP_0$ such that for all $\T\in D$
$$ \sum_{n} \sum_{s\in {\overrightarrow{\T}}^l\cap \OP_n} \left|\langle f_n,\Phi_{s_j}\rangle\right|^2 \geq 2^{2k}|I_\T|$$ and for all sub-trees $\T'\subset \T$
$$ \sum_{n} \sum_{s\in {\overrightarrow{\T'}}^l\cap\OP_n} \left|\langle f_n,\Phi_{s_j}\rangle\right|^2 \leq 2^{2k+2}|I_{\T'}|.$$
For $g$ a function, we similarly define
$$ \widetilde{energy_2}(g,\OP):=\sup_{k\in\Z} \sup_{D} 2^k \left( \sum_{\T \in D} |I_\T| \right)^{1/2},$$
where we take the supremum over all the collections $D$ of strongly $2$-disjoint trees $\T\subset \OP_0$ such that for all $\T\in D$
$$ \sum_{s\in \T} \left|\langle g,\Phi_{s_j}\rangle\right|^2 \geq 2^{2k}|I_\T|$$ and for all sub-trees $\T'\subset \T$
$$ \sum_{s\in \T'} \left|\langle g,\Phi_{s_j}\rangle\right|^2 \leq 2^{2k+2}|I_{\T'}|.$$
\end{df}

\begin{rem} \label{remgg} For the function $g$, the two quantities $size_2$ and $\widetilde{energy_2}$ are bounded by the classical ones (defined for example in \cite{MTT3}.
\end{rem}

\begin{thm} \label{thmsizev} Let $\OP$ be a collection of tri-tiles, $E_j$ be a set of finite measure, $f:=(f_n)$ be a sequence of functions with $\sum_{n} |f_n|^2\leq {\bf 1}_{E_j}$ and $(j,l)$ be two indices with $j\in\{1,3\}$. Then for $p_j>1$
\be{sizev1} \overrightarrow{size_j}^{l}(f,\OP) \lesssim \left(\sup_{s\in \OP} \frac{1}{|I_s|} \int_{E_j}\left(1+ \frac{d(x,I_s)}{|I_s|} \right)^{-100} dx \right)^{1/p_j}. \ee
\end{thm}

\dem The proof is a mixture of the one for Theorem \ref{thm:sizevect1} and Theorem \ref{thm:sizevect3}. That is why we omit details.
We use the following operators (associated to a $k$ tree with $k\neq j$)
$$ Op_n(f_n):=\left(\sum_{s\in {\overrightarrow{\T}}^{l} \cap \OQ_n} \left|\langle f_n,\Phi_{s_j}\rangle \frac{1}{|I_s|^{1/2}}\left(1+\frac{d(x,I_s)}{|I_s|}\right)^{-200}\right|^2 \right)^{1/2}.$$
This vector-valued operator is a modulated Calder\'on-Zygmund operator as the collection ${\overrightarrow{\T}}^{l} \cap \OQ_n$ is a ``classical'' tree (see Lemma \ref{lem:CZ}). Using Theorem 9.5.10 of \cite{Gra} about vector-valued inequality for such operators, we have for all $p>1$~:
$$ \left\| \left( \sum_{n} \left|Op_n(f_n) \right|^2 \right)^{1/2} \right\|_{p} \lesssim \left\| \left( \sum_{n} \left|f_n \right|^2 \right)^{1/2} \right\|_{p}.$$
Then we use previous arguments (explained in details in previously cited Theorems) to conclude the proof.
\findem

\mb For the ``energy'' quantity, we have~:

\begin{thm} \label{thmenerv}
Let $\OP$ be a collection of tri-tile and, $f:=(f_n)_n$ be a sequence of functions and $j\in\{1,3\}$ and $l$ be two indices. Then we have the following estimate for the ``energy'' quantity~:
\be{energy300} \widetilde{{\overrightarrow{energy_j}}^l}(f,\OP) \lesssim \left\|\left(\sum_{n\in\Z} |h_n|^2 \right)^{1/2} \right\|_2. \ee
\end{thm}

\dem This result is included in Theorem \ref{thm:energy3} as for $\T$ a tree, the collection ${\overrightarrow{\T}}^{l}$ is a sub-collection of ${\overrightarrow{\T}}^{1,2}$. So we do not repeat it. 
\findem

\mb After having described the estimates of the different quantities, we have seen in Subsection \ref{subsection2} that the principal result in a consequence of a ``tree estimate'' and a ``good'' algorithm. \\
We begin by the ``tree estimate''~:

\begin{prop} Let $\T$ be a tree of a collection $\OP$ then for every index $l$
$$ \Lambda_{{\overrightarrow{\T}}^{l}}(f_1,f_2,f_3) \lesssim |I_\T| \left[\prod_{i=1,3} {\overrightarrow{size_i}}^l(f_i,\OP) \right] size_2(g,\OP).$$
\end{prop}

\mb The proof is very similar to Proposition \ref{treeesti1}'s. We let the details to the reader.

\mb Now to conclude the proof of Theorem \ref{thm:principal100}, we have just to describe an appropriate algorithm~:

\begin{thm} \label{thm:algo100} Let $j\in\{1,3\}$ and $l\in\{1,2,3\}$ be fixed indices and $f:=(f_n)_n$ be a sequence of functions. Let $\OP$ be a collection of tiles satisfying (\ref{propOQ100}) and such that
 $$ {\overrightarrow{size_j}}^{l}(f_j,\OP) \leq 2^{-n}  \widetilde{{\overrightarrow{energy_j}}^l}(f_j,\OP).$$
Then we may decompose $\OP=\OP^1 \cup \OP^2$ with a collection $\OP^1$ satisfying (\ref{propOQ100}) and
 \be{prop100} {\overrightarrow{size_j}}^{l}(f_j,\OP^1) \leq 2^{-n-1}  \widetilde{{\overrightarrow{energy_j}}^l}(f_j,\OP) \ee
and $\OP^2$ is a collection of vectorial trees $\OP^2=({\overrightarrow{\T_i}}^{1,2})_i$ such that
\be{prop200} \sum_{i} |I_{\T_i}| \lesssim 2^{2n}. \ee
\end{thm}

\mb The proof follows the same reasoning as for Theorem \ref{thm:algo1}. By translation, we use the ``classical'' algorithm in the ``strip of reference'' $\OQ_0$ (which is of rank one). We compute the same arguments as for Theorem \ref{thm:algo1}.

\begin{rem} \label{remg} We have a similar result for the function $g$ with the index $j=2$. In fact, it corresponds to the ``classical'' algorithm in the strip $\OP_0$. So we refer the reader to Proposition 12.2 of \cite{MTT3}.
 \end{rem}

\gb Then Theorem \ref{thm:principal100} comes from the classical analysis. As we have a precise estimate for the ``size quantity'' (see Theorem \ref{thmsizev}) we can prove `` restricted weak type estimates'' for our trilinear form $\Lambda_{\OQ}$ for all the exponents in the described range. We refer the reader to \cite{MTT2} for details. For an easy reference, we just remember the main ideas. 

\gb {\bf Proof of Theorem \ref{thm:principal100}.} \\
The exponents $(p_\beta)_\beta$ and the index $\alpha\in\{1,2,3\}$ are fixed for the proof. Let $E_1,E_2$ and $E_3$ measurable sets of finite measure. First we construct the substantial subset $E_\alpha' \subset E_\alpha$. Denote
$$ U:=\bigcup_{i=1}^3 \left\{ x\in \R,\ M_{HL}({\bf 1}_{E_i})(x)>\eta\frac{|E_i|}{|E_\alpha|}\right\}.$$
By using Hardy-Littlewood Theorem, there exists a numerical constant $\eta$ such that
$$ |U|\leq |E_\alpha|/2.$$
We set also $E_\alpha'=E_\alpha \setminus U$. Now we fix the functions $f_n,g,h_n$ and we shall prove
the inequality (\ref{wtyperestreint}). Using Theorem \ref{thm:algo100}, Remark \ref{remg},  we can use the proof of Theorem \ref{thm:abstrait} to obtain similar results. We fix the index $l=1$ for example. Then for $0<\theta_1,\theta_2,\theta_3$ three real numbers satisfying
$$ \theta_1+\theta_2+\theta_3=1,$$
there exists a constant $C=C(\theta)$ such that for all finite collections $\OQ$ of tri-tiles, all functions $f_1:=(f_{1,n})_n,f_2,f_3:=(f_{3,n})_n$
\begin{align*}
\Lambda_{\OQ}(f_1,f_2,f_3) & \leq C \prod_{j=1,3} \left[\widetilde{{\overrightarrow{energy_j}}^1}(f_j,\OQ)\right]^{1-\theta_j} \left[{\overrightarrow{size_j}}^{1}(f_j,\OQ)\right]^{\theta_j} \\
 & \hspace{2cm} \left[\widetilde{energy_2}(f_2,\OQ)\right]^{1-\theta_2} \left[size_2(f_2,\OQ)\right]^{\theta_2}.
\end{align*}
Then we chose $\theta_\alpha=\frac{2}{p_\alpha}+1 \in(0,1)$ and for $\beta\neq \alpha$, we set $\theta_\beta = \frac{2}{p_\beta}-1 \in(0,1)$. Then it is well-known that the special choice of the substantial subset $E'_\alpha$ with the precise estimates (\ref{sizev1}) and (\ref{energy300}) permits us to obtain the desired estimate
$$\Lambda_{\OQ}(f_1,f_2,f_3) \lesssim \prod_{i=1}^3 |E_i|^{1/p_i}.$$ \findem

\subsection{Other square functions.}

We present here results concerning ``non smooth'' versions of bilinear square functions appeared in Section 4 of \cite{Diestel}.

\mb In the previous estimates, we only have considered a nonsmooth decomposition of the frequency plane by parallel strips. We could be interested by decompositions with other geometric figures, for example parallelograms. Let us explain this example. We consider now two singular angles $\theta_1,\theta_2$, which are supposed to be non degenerate ($\theta_i\in(-\pi/2,\pi/2) \setminus\{0,-\pi/4\}$). We consider now parallelograms defining as follows~: let $a:=(a_n)_n$, $b:=(b_n)_n$, $c:=(c_n)_n$ and $d:=(d_n)_n$ be nondecreasing sequences of numbers satisfying~:
$$  a_n <b_n \leq a_{n+1} \qquad c_n<d_n\leq c_{n+1}.$$
For two integers $n,p$, we construct the following parallelogram~:
$$C_{n,p}:=\left\{ (\xi_1,\xi_2), a_n \leq \xi_2 - \tan(\theta_1)\xi_1\leq b_n \ \textrm{et} \ c_p \leq \xi_2 - \tan(\theta_2)\xi_1\leq d_p\right\}.$$

\begin{figure}[h!]
\centering
\includegraphics{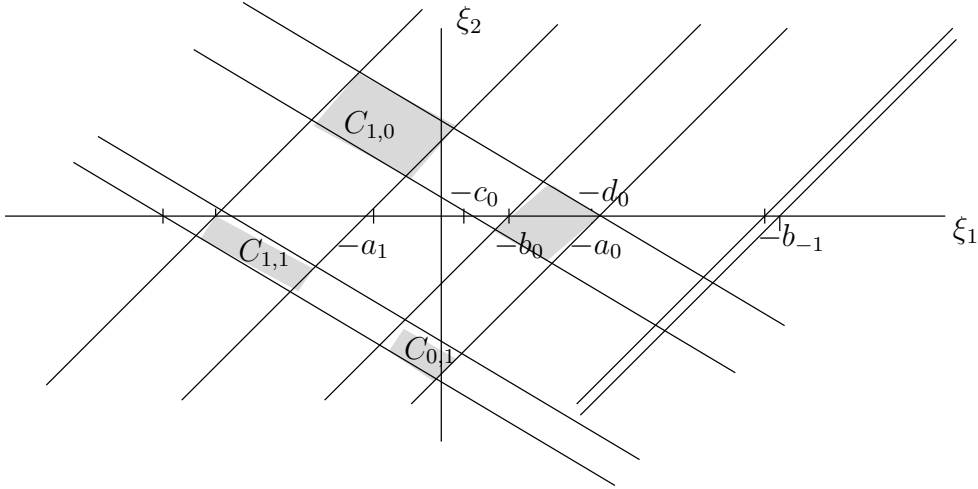}
 \caption{The parallelograms.}
 \label{figure2}
\end{figure}

\mb Then using this new decomposition of the frequency plane, we set $\pi_{C_{n,p}}$ the bilinear multiplier defined by~:
$$\pi_{C_{n,p}}(f,g)(x):= \int_{\R^2} e^{ix(\xi_1+\xi_2)} \widehat{f}(\xi_1) \widehat{g}(\xi_2) {\bf 1}_{C_{n,p}}(\xi) d\xi.$$
We have the following theorem~:

\begin{thm}  \label{thm:principalfin} Let $2<p,q,r'<\infty$ be exponents satisfying
 $$\frac{1}{r}=\frac{1}{p}+\frac{1}{q}.$$
We assume that the sequences $a$ and $b$ satisfy that for all $n\in\Z$
$$ (b_n-a_n)=(b_{n-1}-a_{n-1}) \qquad (a_{n+1}-b_n)=(a_n-{b_{n-1}})$$
and similarly for $c$ and $d$.
Then, there is a constant $C=C(p,q,r)$ such that for all functions $f,g\in\s(\R)$
$$ \left\| \left( \sum_{n,p\in \Z} \left|\pi_{C_{n,p}}(f,g) \right|^2 \right)^{1/2} \right\|_{r} \leq C \|f\|_p \|g\|_q.$$
\end{thm}

\dem We do not write in details the proof of this result. We first decompose each parallelograms with symbols supported in a cone (as explained in Section 4 of \cite{Diestel}). Then we fix a ``parallelogram of reference'' and consider the other ones as translated of this last one. In this set of reference, we can apply the classical algorithm, as the associated collection of tri-tiles is of rank one. We apply the same use of vectorized trees and then similar arguments permit to prove what we want. \findem

\gb We can exactly do the same with other polygons sequences.

\mb To finish, we remember questions, which still stay open. About Theorem \ref{thm:principal}, we do not know if the assumption $2\leq r'$ is necessary or not ? Moreover, it will be very interesting to extend Theorem \ref{thm:principal} for an arbitrary collection of disjoint strips.


\bibliographystyle{plain}
\bibliography{square}

\end{document}